\pdfoutput=1

\documentclass[12pt]{article}

\usepackage{amsfonts,amssymb,amsthm,amsmath,latexsym}
\usepackage{subeqnarray,eqsection,indent,url,cite,bm}
\usepackage{multirow}  
\usepackage{tensor}  
\usepackage{graphicx}
\usepackage{tikz}

\usepackage{color}

\def\red{\textcolor{red}}

\usepackage{framed}  
\renewenvironment{framed}[1][\hsize]
   {\MakeFramed{\hsize#1\advance\hsize-\width \FrameRestore}}%
   {\endMakeFramed}

\date{June 30, 2022 \\[1.5mm] revised December 15, 2022}   

\begin{document}

\title{\vspace*{-1cm}
       A simple algorithm for expanding a power series \\
          as a continued fraction
      }

\author{
     \\
     {\small Alan D.~Sokal}                 \\[2mm]
     {\small\it Department of Mathematics}   \\[-2mm]
     {\small\it University College London}   \\[-2mm]
     {\small\it London WC1E 6BT}             \\[-2mm]
     {\small\it UNITED KINGDOM}              \\[2mm]
     {\small\it Department of Physics}       \\[-2mm]
     {\small\it New York University}         \\[-2mm]
     {\small\it 726 Broadway}                \\[-2mm]
     {\small\it New York, NY 10003}          \\[-2mm]
     {\small\it USA}                         \\[2mm]
     {\small\tt sokal@nyu.edu}               \\[-2mm]
     {\protect\makebox[5in]{\quad}}  
     \\
}

\maketitle
\thispagestyle{empty}   

\begin{abstract}
\noindent
I present and discuss an extremely simple algorithm
for expanding a formal power series as a continued fraction.
This algorithm, which goes back to Euler (1746) and Viscovatov (1805),
deserves to be better known.
I also discuss the connection of this algorithm
with the work of Gauss (1812), Stieltjes (1889), Rogers (1907) and Ramanujan,
and a combinatorial interpretation based on the work of Flajolet (1980).
\end{abstract}


\bigskip
\noindent
{\bf Mathematics Subject Classification (MSC 2010) codes:}
30B70 (Primary); 05A10, 05A15, 05A19 (Secondary).

\clearpage


\theoremstyle{plain} 
\newtheorem{theorem}{Theorem}[section]
\newtheorem{proposition}[theorem]{Proposition}
\newtheorem{lemma}[theorem]{Lemma}
\newtheorem{corollary}[theorem]{Corollary}
\newtheorem{conjecture}[theorem]{Conjecture}
\newtheorem{pseudo-theorem}[theorem]{Pseudo-theorem}

\theoremstyle{definition}
\newtheorem{definition}[theorem]{Definition}
\newtheorem{example}[theorem]{Example}
\newtheorem{question}[theorem]{Question}
\newtheorem{problem}[theorem]{Problem}
\newtheorem{openproblem}[theorem]{Open Problem}

\theoremstyle{remark}
\newtheorem{remark}[theorem]{Remark}

\renewcommand{\theenumi}{\alph{enumi}}
\renewcommand{\labelenumi}{(\theenumi)}
\def\eop{\hbox{\kern1pt\vrule height6pt width4pt
depth1pt\kern1pt}\medskip}
\def\prf{\par\noindent{\bf Proof.\enspace}\rm}
\def\rmk{\par\medskip\noindent{\bf Remark\enspace}\rm}

\newcommand{\textbfit}[1]{\textbf{\textit{#1}}}

\newcommand{\bigdash}{%
\smallskip\begin{center} \rule{5cm}{0.1mm} \end{center}\smallskip}

\newcommand{\be}{\begin{equation}}
\newcommand{\ee}{\end{equation}}
\newcommand{\<}{\langle}
\renewcommand{\>}{\rangle}
\newcommand{\widebar}{\overline}
\def\reff#1{(\protect\ref{#1})}
\def\spose#1{\hbox to 0pt{#1\hss}}
\def\ltapprox{\mathrel{\spose{\lower 3pt\hbox{$\mathchar"218$}}
    \raise 2.0pt\hbox{$\mathchar"13C$}}}
\def\gtapprox{\mathrel{\spose{\lower 3pt\hbox{$\mathchar"218$}}
    \raise 2.0pt\hbox{$\mathchar"13E$}}}
\def\textprime{${}^\prime$}
\def\proof{\par\medskip\noindent{\sc Proof.\ }}
\def\firstproof{\par\medskip\noindent{\sc First Proof.\ }}
\def\secondproof{\par\medskip\noindent{\sc Second Proof.\ }}
\renewcommand{\qed}{ $\square$ \bigskip}
\newcommand{\myendremark}{ $\blacksquare$ \bigskip}
\def\proofof#1{\bigskip\noindent{\sc Proof of #1.\ }}
\def\firstproofof#1{\bigskip\noindent{\sc First Proof of #1.\ }}
\def\secondproofof#1{\bigskip\noindent{\sc Second Proof of #1.\ }}
\def\thirdproofof#1{\bigskip\noindent{\sc Third Proof of #1.\ }}
\def\half{ {1 \over 2} }
\def\third{ {1 \over 3} }
\def\twothird{ {2 \over 3} }
\def\smfrac#1#2{{\textstyle{#1\over #2}}}
\def\smhalf{ {\smfrac{1}{2}} }
\newcommand{\disfrac}[2]{{\displaystyle \frac{#1}{#2}}}
\newcommand{\real}{\mathop{\rm Re}\nolimits}
\renewcommand{\Re}{\mathop{\rm Re}\nolimits}
\newcommand{\imag}{\mathop{\rm Im}\nolimits}
\renewcommand{\Im}{\mathop{\rm Im}\nolimits}
\newcommand{\sgn}{\mathop{\rm sgn}\nolimits}
\newcommand{\tr}{\mathop{\rm tr}\nolimits}
\newcommand{\supp}{\mathop{\rm supp}\nolimits}
\newcommand{\disc}{\mathop{\rm disc}\nolimits}
\newcommand{\diag}{\mathop{\rm diag}\nolimits}
\newcommand{\tridiag}{\mathop{\rm tridiag}\nolimits}
\newcommand{\AZ}{\mathop{\rm AZ}\nolimits}
\newcommand{\perm}{\mathop{\rm perm}\nolimits}
\def\hboxscript#1{ {\hbox{\scriptsize\em #1}} }
\renewcommand{\emptyset}{\varnothing}
\newcommand{\eqdef}{\stackrel{\rm def}{=}}

\newcommand{\restrict}{\upharpoonright}

\newcommand{\compinv}{{\langle -1 \rangle}}   

\newcommand{\scra}{{\mathcal{A}}}
\newcommand{\scrb}{{\mathcal{B}}}
\newcommand{\scrc}{{\mathcal{C}}}
\newcommand{\scrd}{{\mathcal{D}}}
\newcommand{\scre}{{\mathcal{E}}}
\newcommand{\scrf}{{\mathcal{F}}}
\newcommand{\scrg}{{\mathcal{G}}}
\newcommand{\scrh}{{\mathcal{H}}}
\newcommand{\scri}{{\mathcal{I}}}
\newcommand{\scrj}{{\mathcal{J}}}
\newcommand{\scrk}{{\mathcal{K}}}
\newcommand{\scrl}{{\mathcal{L}}}
\newcommand{\scrm}{{\mathcal{M}}}
\newcommand{\scrn}{{\mathcal{N}}}
\newcommand{\scro}{{\mathcal{O}}}
\newcommand{\scrp}{{\mathcal{P}}}
\newcommand{\scrq}{{\mathcal{Q}}}
\newcommand{\scrr}{{\mathcal{R}}}
\newcommand{\scrs}{{\mathcal{S}}}
\newcommand{\scrt}{{\mathcal{T}}}
\newcommand{\scrv}{{\mathcal{V}}}
\newcommand{\scrw}{{\mathcal{W}}}
\newcommand{\scrz}{{\mathcal{Z}}}

\newcommand{\ahat}{{\widehat{a}}}
\newcommand{\Zhat}{{\widehat{Z}}}
\renewcommand{\k}{{\mathbf{k}}}
\newcommand{\n}{{\mathbf{n}}}
\newcommand{\vv}{{\mathbf{v}}}
\newcommand{\bv}{{\mathbf{v}}}
\newcommand{\w}{{\mathbf{w}}}
\newcommand{\x}{{\mathbf{x}}}
\newcommand{\cc}{{\mathbf{c}}}
\newcommand{\zero}{{\mathbf{0}}}
\newcommand{\one}{{\mathbf{1}}}
\newcommand{\bmm}{ {\bf m} }

\newcommand{\C}{{\mathbb C}}
\newcommand{\D}{{\mathbb D}}
\newcommand{\Z}{{\mathbb Z}}
\newcommand{\N}{{\mathbb N}}
\newcommand{\Q}{{\mathbb Q}}
\newcommand{\PP}{{\mathbb P}}
\newcommand{\R}{{\mathbb R}}
\newcommand{\RR}{{\mathbb R}}
\newcommand{\E}{{\mathbb E}}

\newcommand{\Sym}{{\mathfrak{S}}}
\newcommand{\SymB}{{\mathfrak{B}}}

\newcommand{\myle}{\preceq}
\newcommand{\myge}{\succeq}
\newcommand{\mygt}{\succ}

\newcommand{\B}{{\sf B}}
\newcommand{\OB}{{\sf OB}}
\newcommand{\OS}{{\sf OS}}
\newcommand{\OO}{{\sf O}}
\newcommand{\SP}{{\sf SP}}
\newcommand{\OSP}{{\sf OSP}}
\newcommand{\Eu}{{\sf Eu}}
\newcommand{\ERR}{{\sf ERR}}
\newcommand{\sfB}{{\sf B}}
\newcommand{\sfD}{{\sf D}}
\newcommand{\sfE}{{\sf E}}
\newcommand{\sfG}{{\sf G}}
\newcommand{\sfJ}{{\sf J}}
\newcommand{\sfP}{{\sf P}}
\newcommand{\sfQ}{{\sf Q}}
\newcommand{\sfS}{{\sf S}}
\newcommand{\sfT}{{\sf T}}
\newcommand{\sfW}{{\sf W}}
\newcommand{\sfMV}{{\sf MV}}
\newcommand{\AMV}{{\sf AMV}}
\newcommand{\BM}{{\sf BM}}

\newcommand{\emIB}{{\hbox{\em IB}}}
\newcommand{\emIP}{{\hbox{\em IP}}}
\newcommand{\emOB}{{\hbox{\em OB}}}
\newcommand{\emSC}{{\hbox{\em SC}}}

\newcommand{\stat}{{\rm stat}}
\newcommand{\cyc}{{\rm cyc}}
\newcommand{\Asc}{{\rm Asc}}
\newcommand{\asc}{{\rm asc}}
\newcommand{\Des}{{\rm Des}}
\newcommand{\des}{{\rm des}}
\newcommand{\Exc}{{\rm Exc}}
\newcommand{\exc}{{\rm exc}}
\newcommand{\Wex}{{\rm Wex}}
\newcommand{\wex}{{\rm wex}}
\newcommand{\Fix}{{\rm Fix}}
\newcommand{\fix}{{\rm fix}}
\newcommand{\lrmax}{{\rm lrmax}}
\newcommand{\rlmax}{{\rm rlmax}}
\newcommand{\Rec}{{\rm Rec}}
\newcommand{\rec}{{\rm rec}}
\newcommand{\Arec}{{\rm Arec}}
\newcommand{\arec}{{\rm arec}}
\newcommand{\ERec}{{\rm ERec}}
\newcommand{\erec}{{\rm erec}}
\newcommand{\EArec}{{\rm EArec}}
\newcommand{\earec}{{\rm earec}}
\newcommand{\recarec}{{\rm recarec}}
\newcommand{\nonrec}{{\rm nonrec}}
\newcommand{\Cpeak}{{\rm Cpeak}}
\newcommand{\cpeak}{{\rm cpeak}}
\newcommand{\Cval}{{\rm Cval}}
\newcommand{\cval}{{\rm cval}}
\newcommand{\Cdasc}{{\rm Cdasc}}
\newcommand{\cdasc}{{\rm cdasc}}
\newcommand{\Cddes}{{\rm Cddes}}
\newcommand{\cddes}{{\rm cddes}}
\newcommand{\cdrise}{{\rm cdrise}}
\newcommand{\cdfall}{{\rm cdfall}}
\newcommand{\Peak}{{\rm Peak}}
\newcommand{\peak}{{\rm peak}}
\newcommand{\Val}{{\rm Val}}
\newcommand{\val}{{\rm val}}
\newcommand{\Dasc}{{\rm Dasc}}
\newcommand{\dasc}{{\rm dasc}}
\newcommand{\Ddes}{{\rm Ddes}}
\newcommand{\ddes}{{\rm ddes}}
\newcommand{\inv}{{\rm inv}}
\newcommand{\maj}{{\rm maj}}
\newcommand{\rs}{{\rm rs}}
\newcommand{\cross}{{\rm cr}}
\newcommand{\crosshat}{{\widehat{\rm cr}}}
\newcommand{\nest}{{\rm ne}}
\newcommand{\rodd}{{\rm rodd}}
\newcommand{\reven}{{\rm reven}}
\newcommand{\lodd}{{\rm lodd}}
\newcommand{\leven}{{\rm leven}}
\newcommand{\sg}{{\rm sg}}
\newcommand{\bl}{{\rm bl}}
\newcommand{\tran}{{\rm tr}}
\newcommand{\area}{{\rm area}}
\newcommand{\ret}{{\rm ret}}
\newcommand{\peaks}{{\rm peaks}}
\newcommand{\hl}{{\rm hl}}
\newcommand{\sll}{{\rm sll}}
\newcommand{\negg}{{\rm neg}}
\newcommand{\imp}{{\rm imp}}

\newcommand{\ba}{{\bm{a}}}
\newcommand{\bahat}{{\widehat{\bm{a}}}}
\newcommand{\sfa}{{{\sf a}}}
\newcommand{\bb}{{\bm{b}}}
\newcommand{\bc}{{\bm{c}}}
\newcommand{\bff}{{\bm{f}}}
\newcommand{\bg}{{\bm{g}}}
\newcommand{\bll}{{\bm{\ell}}}
\newcommand{\br}{{\bm{r}}}
\newcommand{\bu}{{\bm{u}}}
\newcommand{\bz}{{\bm{z}}}
\newcommand{\bA}{{\bm{A}}}
\newcommand{\bB}{{\bm{B}}}
\newcommand{\bC}{{\bm{C}}}
\newcommand{\bE}{{\bm{E}}}
\newcommand{\bF}{{\bm{F}}}
\newcommand{\bI}{{\bm{I}}}
\newcommand{\bJ}{{\bm{J}}}
\newcommand{\bM}{{\bm{M}}}
\newcommand{\bN}{{\bm{N}}}
\newcommand{\bP}{{\bm{P}}}
\newcommand{\bQ}{{\bm{Q}}}
\newcommand{\bS}{{\bm{S}}}
\newcommand{\bT}{{\bm{T}}}
\newcommand{\bW}{{\bm{W}}}
\newcommand{\bIB}{{\bm{IB}}}
\newcommand{\bOB}{{\bm{OB}}}
\newcommand{\bOS}{{\bm{OS}}}
\newcommand{\bERR}{{\bm{ERR}}}
\newcommand{\bSP}{{\bm{SP}}}
\newcommand{\bMV}{{\bm{MV}}}
\newcommand{\bBM}{{\bm{BM}}}
\newcommand{\balpha}{{\bm{\alpha}}}
\newcommand{\bbeta}{{\bm{\beta}}}
\newcommand{\bgamma}{{\bm{\gamma}}}
\newcommand{\bdelta}{{\bm{\delta}}}
\newcommand{\bomega}{{\bm{\omega}}}
\newcommand{\bone}{{\bm{1}}}
\newcommand{\bzero}{{\bm{0}}}

\newcommand{\Cbar}{{\overline{C}}}
\newcommand{\Dbar}{{\overline{D}}}
\newcommand{\dbar}{{\overline{d}}}
\def\Ctilde{{\widetilde{C}}}
\def\Ftilde{{\widetilde{F}}}
\def\Gtilde{{\widetilde{G}}}
\def\Htilde{{\widetilde{H}}}
\def\Ptilde{{\widetilde{P}}}
\def\Chat{{\widehat{C}}}
\def\ctilde{{\widetilde{c}}}
\def\zbar{{\overline{Z}}}
\def\pitilde{{\widetilde{\pi}}}

\newcommand{\sech}{{\rm sech}}

%
%
\newcommand{\sn}{{\rm sn}}
\newcommand{\cn}{{\rm cn}}
\newcommand{\dn}{{\rm dn}}
\newcommand{\sm}{{\rm sm}}
\newcommand{\cm}{{\rm cm}}

%
%
\newcommand{\zfz}{ {{}_0 \! F_0} }
\newcommand{\oft}{ {{}_1 \! F_2} }

%
%
\newcommand{\FHyper}[2]{ {\tensor[_{#1 \!}]{F}{_{#2}}\!} }
\newcommand{\FHYPER}[5]{ {\FHyper{#1}{#2} \!\biggl(
   \!\!\begin{array}{c} #3 \\[1mm] #4 \end{array}\! \bigg|\, #5 \! \biggr)} }
\newcommand{\tfo}{ {\FHyper{2}{1}} }
\newcommand{\tfz}{ {\FHyper{2}{0}} }
\newcommand{\ofo}{ {\FHyper{1}{1}} }
\newcommand{\ofz}{ {\FHyper{1}{0}} }
\newcommand{\zfo}{ {\FHyper{0}{1}} }
\newcommand{\FHYPERbottomzero}[3]{ {\FHyper{#1}{0} \!\biggl(
   \!\!\begin{array}{c} #2 \\[1mm] \hbox{---} \end{array}\! \bigg|\, #3 \! \biggr)} }
\newcommand{\FHYPERtopzero}[3]{ {\FHyper{0}{#1} \!\biggl(
   \!\!\begin{array}{c} \hbox{---} \\[1mm] #2 \end{array}\! \bigg|\, #3 \! \biggr)} }

\newcommand{\phiHyper}[2]{ {\tensor[_{#1}]{\phi}{_{#2}}} }
\newcommand{\psiHyper}[2]{ {\tensor[_{#1}]{\psi}{_{#2}}} }
\newcommand{\PhiHyper}[2]{ {\tensor[_{#1}]{\Phi}{_{#2}}} }
\newcommand{\PsiHyper}[2]{ {\tensor[_{#1}]{\Psi}{_{#2}}} }
\newcommand{\phiHYPER}[6]{ {\phiHyper{#1}{#2} \!\left(
   \!\!\begin{array}{c} #3 \\ #4 \end{array}\! ;\, #5, \, #6 \! \right)\!} }
\newcommand{\psiHYPER}[6]{ {\psiHyper{#1}{#2} \!\left(
   \!\!\begin{array}{c} #3 \\ #4 \end{array}\! ;\, #5, \, #6 \! \right)} }
\newcommand{\PhiHYPER}[5]{ {\PhiHyper{#1}{#2} \!\left(
   \!\!\begin{array}{c} #3 \\ #4 \end{array}\! ;\, #5 \! \right)\!} }
\newcommand{\PsiHYPER}[5]{ {\PsiHyper{#1}{#2} \!\left(
   \!\!\begin{array}{c} #3 \\ #4 \end{array}\! ;\, #5 \! \right)\!} }
\newcommand{\zerophizero}{ {\phiHyper{0}{0}} }
\newcommand{\ophizero}{ {\phiHyper{1}{0}} }
\newcommand{\zphio}{ {\phiHyper{0}{1}} }
\newcommand{\ophio}{ {\phiHyper{1}{1}} }
\newcommand{\tphio}{ {\phiHyper{2}{1}} }
\newcommand{\tphiz}{ {\phiHyper{2}{0}} }
\newcommand{\tPhio}{ {\PhiHyper{2}{1}} }
\newcommand{\opsio}{ {\psiHyper{1}{1}} }

%
%
\newcommand{\stirlingsubset}[2]{\genfrac{\{}{\}}{0pt}{}{#1}{#2}}
\newcommand{\stirlingcycleold}[2]{\genfrac{[}{]}{0pt}{}{#1}{#2}}
\newcommand{\stirlingcycle}[2]{\left[\! \stirlingcycleold{#1}{#2} \!\right]}
\newcommand{\assocstirlingsubset}[3]{{\genfrac{\{}{\}}{0pt}{}{#1}{#2}}_{\! \ge #3}}
\newcommand{\genstirlingsubset}[4]{{\genfrac{\{}{\}}{0pt}{}{#1}{#2}}_{\! #3,#4}}
\newcommand{\euler}[2]{\genfrac{\langle}{\rangle}{0pt}{}{#1}{#2}}
\newcommand{\eulergen}[3]{{\genfrac{\langle}{\rangle}{0pt}{}{#1}{#2}}_{\! #3}}
\newcommand{\eulersecond}[2]{\left\langle\!\! \euler{#1}{#2} \!\!\right\rangle}
\newcommand{\eulersecondgen}[3]{{\left\langle\!\! \euler{#1}{#2} \!\!\right\rangle}_{\! #3}}
\newcommand{\binomvert}[2]{\genfrac{\vert}{\vert}{0pt}{}{#1}{#2}}
\newcommand{\binomsquare}[2]{\genfrac{[}{]}{0pt}{}{#1}{#2}}
\newcommand{\qbinom}[3]{\genfrac{(}{)}{0pt}{}{#1}{#2}_{\! #3}}


\newenvironment{sarray}{
             \textfont0=\scriptfont0
             \scriptfont0=\scriptscriptfont0
             \textfont1=\scriptfont1
             \scriptfont1=\scriptscriptfont1
             \textfont2=\scriptfont2
             \scriptfont2=\scriptscriptfont2
             \textfont3=\scriptfont3
             \scriptfont3=\scriptscriptfont3
           \renewcommand{\arraystretch}{0.7}
           \begin{array}{l}}{\end{array}}

\newenvironment{scarray}{
             \textfont0=\scriptfont0
             \scriptfont0=\scriptscriptfont0
             \textfont1=\scriptfont1
             \scriptfont1=\scriptscriptfont1
             \textfont2=\scriptfont2
             \scriptfont2=\scriptscriptfont2
             \textfont3=\scriptfont3
             \scriptfont3=\scriptscriptfont3
           \renewcommand{\arraystretch}{0.7}
           \begin{array}{c}}{\end{array}}



\newcommand{\chapquote}[2]{
  \begin{list}{}{
     \setlength{\leftmargin}{0.35\textwidth}
     \setlength{\rightmargin}{0.05\textwidth}
     \setlength{\labelwidth}{0pt}
  }
  \item
  \footnotesize
  #1
  \par\raggedleft
  {\sl --- #2}
  \end{list}
  \normalsize\smallskip
}

\newcommand{\chapquotenoskip}[2]{
  \begin{list}{}{
     \setlength{\leftmargin}{0.35\textwidth}
     \setlength{\rightmargin}{0.05\textwidth}
     \setlength{\labelwidth}{0pt}
  }
  \item
  \footnotesize
  #1
  \par\raggedleft
  {\sl --- #2}
  \end{list}
  \normalsize
}

\hyphenation{Stieltjes Flajolet Gaurav Mathematica}

\clearpage

\chapquotenoskip{
   Surely the story unfolded here emphasizes how valuable it~is
   to study and understand the central ideas behind major pieces
   of mathematics produced by giants like Euler.
}{
   George Andrews \cite[p.~284]{Andrews_83}
}

\section{Introduction}

The expansion of power series into continued fractions
goes back nearly 300 years.
Euler \cite{Euler_1760} showed circa 1746 that\footnote{
   The paper \cite{Euler_1760},
   which is E247 in Enestr\"om's \cite{Enestrom_13} catalogue,
   was probably written circa 1746;
   it~was presented to the St.~Petersburg Academy in 1753
   and published in 1760.
   See also \cite{Barbeau_76,Barbeau_79,Varadarajan_07}
   for some commentary on the analytic aspects of this paper.
}
\be
   \sum_{n=0}^\infty n! \: t^n
   \;=\;
   \cfrac{1}{1 - \cfrac{1t}{1 - \cfrac{1t}{1 - \cfrac{2t}{1- \cfrac{2t}{1- \cfrac{3t}{1 - \cfrac{3t}{1 - \cdots}}}}}}}
 \label{eq.nfact.contfrac}
\ee
and more generally that
\be
   \sum_{n=0}^\infty a(a+1)(a+2) \cdots (a+n-1) \: t^n
   \;=\;
   \cfrac{1}{1 - \cfrac{at}{1 - \cfrac{1t}{1 - \cfrac{(a+1)t}{1- \cfrac{2t}{1-\cfrac{(a+2)t}{1 - \cfrac{3t}{1- \cdots}}}}}}}
   \;\,.
 \label{eq.xupperfact.contfrac}
\ee
Lambert \cite{Lambert_1768} showed circa 1761 that\footnote{
   Several sources
   (e.g.~\cite{Laczkovich_97,Wallisser_00}
         \cite[p.~110]{Brezinski_91}
         \cite[p.~327]{Legendre_1800})   
   date Lambert's proof to 1761,
   though I~am not sure what is the evidence for this.
   Lambert's paper was read to the Royal Prussian Academy of Sciences
   in 1767, and published in 1768.
   See \cite{Laczkovich_97,Wallisser_00} for analyses of
   Lambert's remarkable work.
}
\be
   {\tan t \over t}
   \;=\;
   \cfrac{1}{1 - \cfrac{{1 \over 1 \cdot 3}t^2}{1 - \cfrac{{1 \over 3 \cdot 5}t^2}{1 - \cfrac{{1 \over 5 \cdot 7}t^2}{1- \cfrac{{1 \over 7 \cdot 9}t^2}{1-\cdots}}}}}
 \label{eq.lambert}
\ee
and used it to prove the irrationality of $\pi$
\cite{Laczkovich_97,Wallisser_00}.\footnote{
   In fact, as noted by Brezinski \cite[p.~110]{Brezinski_91},
   a formula equivalent to \reff{eq.lambert} appears already
   in Euler's first paper on continued fractions \cite{Euler_1744}:
   see top p.~321 in the English translation.
   The paper \cite{Euler_1744},
   which is E71 in Enestr\"om's \cite{Enestrom_13} catalogue,
   was presented to the St.~Petersburg Academy in 1737
   and published in 1744.
}
Many similar expansions were discovered in the eighteenth
and nineteenth centuries:
most notably, Gauss \cite{Gauss_1813} found in 1812
a continued-fraction expansion for the ratio
of two contiguous hypergeometric functions $\tfo$,
from which many previously obtained expansions
can be deduced by specialization or taking limits \cite[Chapter~XVIII]{Wall_48}.
A detailed history of continued fractions can be found
in the fascinating book of Brezinski \cite{Brezinski_91}.

Let us stress that this subject has two facets: algebraic and analytic.
The algebraic theory treats both sides of identities like
\reff{eq.nfact.contfrac}--\reff{eq.lambert}
as formal power series in the indeterminate $t$;
convergence plays no role.\footnote{
   See \cite{Niven_69}, \cite[Chapter~1]{Henrici_74} or
   \cite[Chapter~2]{Wilf_94}
   for an introduction to formal power series;
   and see \cite[Section~IV.4]{Bourbaki_90} for a more complete treatment.
}
Thus, \reff{eq.nfact.contfrac} is a perfectly meaningful (and true!)\ identity
for formal power series, despite the fact that the left-hand side
has zero radius of convergence.
By contrast, the analytic theory seeks to understand the regions
of the complex $t$-plane in which the left or right sides of the identity
are well-defined, to study whether they are equal there,
and to investigate possible analytic continuations.
In this paper we shall be concerned solely with the algebraic aspect;
indeed, the coefficients in our formulae need not be complex numbers,
but may lie in an arbitrary field $F$.

The central goal of this paper is to present and discuss
an extremely simple algorithm for expanding a formal power series
$f(t) = \sum\limits_{n=0}^\infty a_n t^n$ with $a_0 \neq 0$
as a continued fraction of the form
\be
   f(t)  \;=\;
   \cfrac{\alpha_0}{1 - \cfrac{\alpha_1 t^{p_1}}{1 - \cfrac{\alpha_2 t^{p_2}}{1 -   \cdots}}}
   \label{def.Stype}
\ee
with integer powers $p_i \ge 1$,
or more generally as
\be
   f(t)  \;=\;
   \cfrac{\alpha_0}{1 - \sum\limits_{j=1}^{M_1} \delta_1^{(j)} t^j - \cfrac{\alpha_1 t^{p_1}}{1 - \sum\limits_{j=1}^{M_2} \delta_2^{(j)} t^j - \cfrac{\alpha_2 t^{p_2}}{1 - \cdots}}}
   \label{def.Jtype}
\ee
with integers $M_i \ge 0$ and $p_i \ge M_i +1$.
Most generally, we will consider continued fractions of the form
\be
   f(t)  \;=\;
   \cfrac{A_0(t)}{1 - \Delta_1(t) - \cfrac{A_1(t)}{1 - \Delta_2(t) - \cfrac{A_2(t)}{1 - \cdots}}}
   \label{def.cfgen}
\ee
where $A_0(t)$ is a formal power series with nonzero constant term,
and $\Delta_k(t)$ and $A_k(t)$ for $k \ge 1$
are formal power series with zero constant term.

In the classical literature on continued fractions
\cite{Perron,Wall_48,Khovanskii_63,Jones_80,Lorentzen_92,Cuyt_08},
\reff{def.Stype} is called a (general) C-fraction \cite{Leighton_39};
with $p_1 = p_2 = \ldots = 1$ it is called a regular C-fraction;
and \reff{def.Jtype} with
$M_1 = M_2 = \ldots =1$ and $p_1 = p_2 = \ldots = 2$
is called an associated continued fraction.
In the recent combinatorial literature on continued fractions
\cite{Flajolet_80,Viennot_83},
the regular C-fraction is called a Stieltjes-type continued fraction
(or S-fraction),
and the associated continued fraction is called a Jacobi-type continued fraction
(or J-fraction).\footnote{
   In the classical literature on continued fractions,
   the terms ``S-fraction'' and ``J-fraction'' refer to closely related
   but {\em different}\/ objects.
}

Although the algorithm to be presented here is more than two centuries old,
it does not seem to be very well known,
or its simplicity adequately appreciated.
In the special case \reff{eq.nfact.contfrac}
it goes back to Euler in 1746 \cite[section~21]{Euler_1760},
as will be explained in Section~\ref{sec.examples.1} below.
While reading \cite{Euler_1760} I~realized that Euler's method
is in fact a completely general algorithm, applicable to arbitrary
power series (perhaps Euler himself already knew this).
Only later did I learn that a substantially equivalent algorithm
was proposed by Viscovatov \cite{Viscovatov_1806} in 1805
and presented in modern notation in the book of
Khovanskii \cite[pp.~27--31]{Khovanskii_63}.\footnote{
   Some post-Khovanskii books on continued fractions
   also discuss the Viscovatov algorithm
   \cite[pp.~2, 16--17, 89--90]{Cuyt_87}
   \cite[pp.~259--265]{Lorentzen_92}
   \cite[pp.~133--141]{Baker_96}
   \cite[pp.~20--21, 112--113, 118--119]{Cuyt_08},
   but in my opinion they do not sufficiently stress
   its simplicity and importance.
   Viscovatov's work is also discussed briefly in
   Brezinski's historical survey \cite[p.~190]{Brezinski_91}.
   Perron, in his classic monograph \cite{Perron}, mentions in a footnote
   the ``useful recursive formula'' of Viscovatov
   \cite[1st ed., p.~304; 2nd ed., p.~304; 3rd ed., vol.~2, p.~120]{Perron},
   but without further explanation.
   See also the Historical Remark
   at \reff{eq.khovanskii.1}/\reff{eq.khovanskii.2} below.
}
I~therefore refer to it as the \textbfit{Euler--Viscovatov algorithm}.
This algorithm was rediscovered several times in the mid-twentieth century
\cite{Gordon_68} \cite{Watson_73} \cite{ODonohoe_74,Murphy_75,Murphy_77}
and probably many times earlier as well.
I~would be very grateful to any readers
who could point out additional relevant references.

Many other algorithms for expanding a power series as a continued fraction
are known, notably the quotient-difference algorithm
\cite[Section~7.1.2]{Jones_80}.
The key advantage of the algorithm presented here is that it avoids
all nonlinear operations on power series
(such as multiplication or division).

But the Euler--Viscovatov algorithm is more than just an algorithm
for computing continued fractions;
suitably reinterpreted,
it becomes a method for {\em proving}\/ continued fractions.
Since this method
was employed implicitly by Euler \cite[section~21]{Euler_1760}
for proving \reff{eq.nfact.contfrac}
and explicitly by Gauss \cite[sections~12--14]{Gauss_1813}
for proving his continued fraction for ratios of contiguous $\tfo$,
I~shall call it the \textbfit{Euler--Gauss recurrence method},
and I~will illustrate it with a variety of examples.

Unless stated otherwise, I~shall assume that the coefficients
$a_i$, $\alpha_i$ and $\delta_i^{(j)}$ belong to a field $F$.
Later I~shall make some brief remarks about what happens
when the coefficients lie instead in a
commutative ring-with-identity-element $R$.

\section{Expansion as a C-fraction}   \label{sec.Stype}

To each continued fraction of the form \reff{def.Stype}
there manifestly corresponds a unique formal power series
$f(t) = \sum_{n=0}^\infty a_n t^n$;
and clearly $\alpha_0 = 0$ if and only if $f(t)$ is identically zero.
Since we are always assuming that $a_0 \neq 0$,
it follows that $\alpha_0 = a_0 \neq 0$.

We say that a continued fraction of the form \reff{def.Stype}
with $\alpha_0 \neq 0$ is
\textbfit{terminating of length $\bm{k}$} ($k \ge 0$)
if $\alpha_1,\ldots,\alpha_k \neq 0$ and $\alpha_{k+1} = 0$;
we say that it is \textbfit{nonterminating}
if all the $\alpha_i$ are nonzero.
Two continued fractions of the form \reff{def.Stype}
will be called \textbfit{equivalent}
if they are both terminating of the same length $k$
and they have the same values for
$\alpha_1,\ldots,\alpha_k$ and $p_1,\ldots,p_k$
(and of course for $\alpha_{k+1} = 0$);
they then correspond to the same power series $f(t)$,
irrespective of the values of $\alpha_{k+2},\alpha_{k+3},\ldots$
and $p_{k+1},p_{k+2},\ldots\,$, which play no role whatsoever.

We shall use the notation $[t^m] \, g(t)$ to denote the coefficient of $t^m$
in the formal power series $g(t)$.

Given a continued fraction of the form \reff{def.Stype},
let us define for $k \ge 0$
\be
   f_k(t)  \;=\;
   \cfrac{1}{1 - \cfrac{\alpha_{k+1} t^{p_{k+1}}}{1 - \cfrac{\alpha_{k+2} t^{p_{k+2}}}{1 - \cdots}}}
   \;\,;
 \label{def.fk}
\ee
of course these are formal power series with constant term 1.
We thus have $f(t) = \alpha_0 f_0(t)$ and the recurrence
\be
   f_k(t)  \;=\;  {1 \over 1 \,-\, \alpha_{k+1} t^{p_{k+1}} f_{k+1}(t)}
   \qquad\hbox{for $k \ge 0$}  \;.
 \label{eq.recurrence.fk}
\ee
Given $f(t)$, we can reconstruct
$(\alpha_k)_{k \ge 0}$, $(p_k)_{k \ge 1}$ and $(f_k)_{k \ge 0}$
by the following obvious algorithm:

\begin{framed}[0.85\textwidth]
\noindent
{\bf Primitive algorithm.}

\medskip
\noindent
1.  Set $\alpha_0 = a_0 = [t^0] \, f(t)$ and
      $f_0(t) = \alpha_0^{-1} f(t)$.

\medskip
\noindent
2.  For $k=1,2,3,\ldots$, do:
\begin{itemize}
   \item[(a)]  If $f_{k-1}(t) = 1$, set $\alpha_k = 0$ and terminate.
      [Then $\alpha_{k+1},\alpha_{k+2},\ldots$ and $p_{k},p_{k+1},\ldots$
       can be given completely arbitrary values.]
   \item[(b)]  If $f_{k-1}(t) \neq 1$, let $p_k$ be the smallest index $n \ge 1$
       such that $[t^n] \, f_{k-1}(t) \neq 0$;
       set $\alpha_k = [t^{p_k}] \, f_{k-1}(t)$;
       and set
\be
   f_k(t)  \;=\;  \alpha_k^{-1} t^{-p_k} \biggl( 1 \,-\, {1 \over f_{k-1}(t)}
                                         \biggr)
   \;.
 \label{eq.alg.fk}
\ee
\end{itemize}
\vspace*{-5mm}
\end{framed}

If this algorithm terminates, then obviously $f$ is a rational function.
Conversely, if $f$ is a rational function,
then it is not difficult to show,
by looking at the degrees of numerator and denominator,
that the algorithm must terminate.
(I~will give the details of this argument a bit later.)
The algorithm therefore proves:

\begin{proposition}
{$\!\!\!$ \bf (Leighton and Scott \protect\cite{Leighton_39})\ }
   \label{prop.leighton}
Let $f(t) = \sum_{n=0}^\infty a_n t^n$
be a formal power series with coefficients in a field $F$, with $a_0 \neq 0$.
Then $f(t)$ can be represented by a continued fraction of the form
\reff{def.Stype}, which is unique modulo equivalence.
This continued fraction is terminating if and only if
$f(t)$ represents a rational function.
\end{proposition}

The disadvantage of the foregoing algorithm
is that it requires division of power series in the step \reff{eq.alg.fk}.
To eliminate this, let us define
\be
   g_k(t)  \;=\;  \prod_{i=0}^k f_i(t) \qquad\hbox{for $k \ge -1$}  \;;
 \label{def.gk}
\ee
these are formal power series with constant term 1, which satisfy
$g_{-1}(t) = 1$ and
\be
   f_k(t)  \;=\;  {g_k(t) \over g_{k-1}(t)}  \qquad\hbox{for $k \ge 0$}  \;.
 \label{eq.fk.gk}
\ee
Then the nonlinear two-term recurrence \reff{eq.recurrence.fk} for the $(f_k)$
becomes the {\em linear}\/ three-term recurrence
\be
   g_k(t) - g_{k-1}(t)  \;=\; \alpha_{k+1} t^{p_{k+1}} g_{k+1}(t)
 \label{eq.recurrence.gk.0}
\ee
for the $(g_k)$.
Rewriting the algorithm in terms of $(g_k)_{k \ge -1}$, we have:

\begin{framed}[0.88\textwidth]
\noindent
{\bf Refined algorithm.}

\medskip
\noindent
1. Set $g_{-1}(t) = 1$, $\alpha_0 = a_0 = [t^0] \, f(t)$ and
      $g_0(t) = \alpha_0^{-1} f(t)$.

\medskip
\noindent
2. For $k=1,2,3,\ldots$, do:
\begin{itemize}
   \item[(a)]  If $g_{k-1}(t) = g_{k-2}(t)$, set $\alpha_k = 0$ and terminate.
   \item[(b)]  If $g_{k-1}(t) \neq g_{k-2}(t)$, let $p_k$ be the smallest index
       $n$ such that $[t^n] \, g_{k-1}(t) \neq [t^n] \, g_{k-2}(t)$;
       set $\alpha_k = [t^{p_k}] \, \bigl(g_{k-1}(t) - g_{k-2}(t) \bigr)$;
       and set
\be
   g_k(t)  \;=\;  \alpha_k^{-1} t^{-p_k} \bigl(g_{k-1}(t) - g_{k-2}(t) \bigr)
   \;.
 \label{eq.alg.gk}
\ee
\end{itemize}
\vspace*{-5mm}
\end{framed}

\noindent
This algorithm requires only {\em linear}\/ operations on power series
(together, of course, with a nonlinear operation in the field $F$,
 namely, division by $\alpha_k$).

Let us also observe that it is not mandatory to take $g_{-1} = 1$.
In fact, we can let $g_{-1}$ be {\em any}\/ formal power series
with constant term 1, and replace \reff{def.gk} by
\be
   g_k(t)  \;=\;  g_{-1}(t) \prod_{i=0}^k f_i(t) \qquad\hbox{for $k \ge 0$}  \;;
 \label{def.gk.generalized}
\ee
then the key relation \reff{eq.fk.gk} still holds.
The algorithm becomes:
\begin{framed}[0.85\textwidth]
\noindent
{\bf Refined algorithm, generalized version.}

\medskip
\noindent
1.  Choose any formal power series $g_{-1}(t)$ with constant term~1;
      \linebreak
      then set $\alpha_0 = a_0 = [t^0] \, f(t)$ and
      $g_0(t) = \alpha_0^{-1} g_{-1}(t) f(t)$.

\medskip
\noindent
2.  As before.
\end{framed}
\noindent
This generalization is especially useful in case $f(t)$ happens to be
given to us as an explicit fraction;  then we can (if we wish) choose $g_{-1}$
to be the denominator.

In particular, suppose that $f = P/Q$ is a rational function
normalized to $Q(0) = 1$,
and that we choose $g_{-1} = Q$ and $g_0 = P/P(0)$.
Then all the $g_k$ are polynomials, and we have
\be
   \deg g_k
   \;\le\;
   \max(\deg g_{k-1}, \deg g_{k-2}) - p_k
   \;\le\;
   \max(\deg g_{k-1}, \deg g_{k-2}) - 1  \;.
\ee
It follows by induction that
\be
   \deg g_k  \;\le\;  d - \lceil k/2 \rceil
 \label{eq.deg.g_k}
\ee
where $d = \max(\deg P,\deg Q)$ is the degree of $f$.
Hence the algorithm (in any of its versions)
must terminate no later than step $k=2d$.
This completes the proof of Proposition~\ref{prop.leighton}.

Of course, the foregoing algorithms, interpreted literally,
require manipulations on power series with infinitely many terms.
Sometimes this can be done by hand
(as we shall see in Sections~\ref{sec.examples.1}--\ref{sec.examples.3})
or by a sufficiently powerful symbolic-algebra package,
if explicit formulae for the series coefficients are available.
But in many cases we are given the initial series $f(t)$
only through some order $t^N$, and we want to find a continued fraction
of the form \reff{def.Stype} that represents $f(t)$
at least through this order.
This can be done as follows:
We start by writing $f(t) = \sum_{n=0}^N a_n t^n + O(t^{N+1})$
and then carry out the algorithm (in any version)
where each $f_k$ or $g_k$ is written as a finite sum
{\em plus an explicit error term $O(t^{N_k+1})$}\/.\footnote{
   This is how {\sc Mathematica} automatically handles
   {\tt SeriesData} objects,
   and how {\sc Maple} handles the {\tt series} data structure.
}
Clearly $N_k = N - \sum_{i=1}^k p_i$.
The algorithm terminates when
$f_{k-1}(t) = 1 + O(t^{N_{k-1}+1})$ or
$g_{k-1}(t) - g_{k-2}(t) = O(t^{N_{k-1}+1})$.
In terms of the coefficients $g_{k,n}$
in $g_k(t) = \sum_{n=0}^\infty g_{k,n} t^n$ (where $g_{k,0} = 1$),
the refined algorithm (in the generalized version) is as follows:
%
%
\begin{framed}[0.9\textwidth]
\noindent
{\bf Refined algorithm, finite-$\bm{N}$ version.}

\bigskip
\noindent
INPUT: Coefficients $g_{k,n}$ for $k=-1,0$ and $0 \le n \le N$,
         where $g_{-1,0} = g_{0,0} = 1$.

\bigskip
\noindent
1. Set $N_0 = N$.

\medskip
\noindent
2. For $k=1,2,3,\ldots$, do:
\begin{itemize}
   \item[(a)]  If $g_{k-1,n} = g_{k-2,n}$ for $0 \le n \le N_{k-1}$,
       set $\alpha_k = 0$ and terminate.
   \item[(b)]  Otherwise, let $p_k$ be the smallest index $n$ ($\le N_{k-1}$)
       such that $g_{k-1,n} \neq g_{k-2,n}$;
       set $\alpha_k = g_{k-1,p_k} - g_{k-2,p_k}$;
       set $N_k = N_{k-1} - p_k$;
       and set
\be
   g_{k,n}  \;=\;  \alpha_k^{-1} (g_{k-1,n+p_k}- g_{k-2,n+p_k})
   \quad \hbox{for $0 \le n \le N_k$}
   \;.
 \label{eq.alg.gkn}
\ee
\end{itemize}
\vspace*{-5mm}
\end{framed}
\noindent
It is easy to see that this algorithm requires
$O(N^2)$ field operations to find a continued fraction
that represents $f(t)$ through order $t^N$.
Note also that if it is subsequently desired to extend the computation
to larger $N$,
one can return to $k=0$ and compute the new coefficients $g_{k,n}$
using \reff{eq.alg.gkn}, without needing to revisit the old ones;
this is a consequence of the method's linearity.

\bigskip

{\bf Historical remark.}
While preparing this article I learned that the
``refined algorithm'' is essentially equivalent
(when $p_1 = p_2 = \ldots = 1$)
to a method presented by Viscovatov \cite[p.~228]{Viscovatov_1806} in 1805.
In terms of Khovanskii's \cite[pp.~27--31]{Khovanskii_63}
quantities $\alpha_{m,n}$,
it suffices to define
\be
   g_m(t)  \;=\; \sum_{n=0}^\infty {\alpha_{m,n} \over \alpha_{m,0}} \, t^n
 \label{eq.khovanskii.1}
\ee
and
\be
   \alpha_m  \;=\;  - \, {\alpha_{m,0} \over \alpha_{m-1,0} \, \alpha_{m-2,0}}
   \;;
 \label{eq.khovanskii.2}
\ee
then Khovanskii's recurrence
$\alpha_{m,n} =
 \alpha_{m-1,0} \alpha_{m-2,n+1} - \alpha_{m-2,0} \alpha_{m-1,n+1}$
\cite[p.~28]{Khovanskii_63}
is equivalent to our \reff{eq.alg.gkn} specialized to $p_k = 1$.
See also
\cite[p.~547, eqns.~(12.6-26) and (12.6-27)]{Henrici_77},
\cite[p.~17]{Cuyt_87}
and \cite[p.~20, eq.~(1.7.7) and p.~112, eq.~(6.1.12c)]{Cuyt_08}.
This same recurrence was independently discovered in the mid-twentieth century
by Gordon \cite[Appendix~A]{Gordon_68},
who named it the ``product-difference algorithm'';
by P.J.S.~Watson \cite[p.~94]{Watson_73};
and by O'Donohoe \cite{ODonohoe_74,Murphy_75,Murphy_77},
who named it the ``corresponding sequence (CS) algorithm''.
The presentation in \cite[Chapter~3]{ODonohoe_74} \cite{Murphy_77}
is particularly clear.

It should be mentioned, however, that some of the modern works
that refer to the ``Viscovatov algorithm'' fail to distinguish
clearly between the primitive algorithm \reff{eq.alg.fk}
and the refined algorithm \reff{eq.alg.gk}/\reff{eq.alg.gkn}.
However, the modern authors should not be blamed:
Viscovatov \cite{Viscovatov_1806} himself fails to make this distinction clear.
As Khovanskii \cite[p.~27]{Khovanskii_63} modestly says,
``This procedure was {\em in principle}\/ [emphasis mine]
proposed by V.~Viskovatoff;
we have merely developed a more convenient notation for this method
of calculation.''

See also \cite{history} for fascinating information
concerning the life of Vasili\u{\i} Ivanovich Viscovatov (1780--1812).

A very similar algorithm was presented by
Christoph (or Christian) Friedrich Kausler (1760--1825)
in 1802 \cite{Kausler_1806}
(see also \cite[pp.~112 ff.]{Kausler_1803}),
but the precise relation between the two algorithms is not clear to me.
\myendremark

\bigskip

We can also run this algorithm in reverse.
Suppose that we have a sequence $(g_k)_{k \ge -1}$
of formal power series with constant term 1,
which satisfy a recurrence of the form
\be
   g_k(t) - g_{k-1}(t)  \;=\;  \alpha_{k+1} t^{p_{k+1}} g_{k+1}(t)
   \qquad\hbox{for $k \ge 0$}  \;.
 \label{eq.gk.recurrence.general}
\ee
(We need not assume that $g_{-1} = 1$.)
Then the series $(f_k)_{k \ge 0}$ defined by $f_k = g_k/g_{k-1}$
satisfy the recurrence \reff{eq.recurrence.fk};
and iterating this recurrence, we see that they are given
by the continued fractions \reff{def.fk}.
This method for proving continued fractions
was employed implicitly by Euler \cite[section~21]{Euler_1760}
for proving \reff{eq.nfact.contfrac}
--- as we shall see in the next section ---
and explicitly by Gauss \cite[sections~12--14]{Gauss_1813}
for proving his continued fraction for ratios of contiguous $\tfo$.
We therefore call it the \textbfit{Euler--Gauss recurrence method}.

\bigskip

Suppose, finally, that the coefficients of $f(t)$ lie in a
commutative ring-with-identity-element $R$, not necessarily a field.
There are two cases:

(a) If $R$ is an integral domain (i.e.\ has no divisors of zero),
then we can carry out the preceding algorithm (in any version)
in the field of fractions $F(R)$,
yielding coefficients $(\alpha_k)_{k \ge 0}$ that lie in $F(R)$
and are unique modulo equivalence.
In some cases these coefficients may lie in $R$, in other cases not.
If $(\alpha_k)_{k \ge 0}$ do lie in $R$, then so will all the coefficients
of the series $(f_k)_{k \ge 0}$ and $(g_k)_{k \ge 0}$.\footnote{
   I am assuming here that the coefficients of the chosen $g_{-1}$
   lie in $R$.
}
In this case the algorithm can be carried out entirely in $R$;
it requires divisions $a/b=c$, but only in cases where $c$
lies in $R$ (and $c$ is of course unique because $R$ has no divisors of zero).

(b) If, by contrast, $R$ has divisors of zero, then the expansion as a
continued fraction can be highly nonunique.
For instance, in $R = \Z_4$, the series $f(t) = 1 + 2t$
is represented in the form \reff{def.Stype} with $p_1 = p_2 = \ldots = 1$
for {\em any}\/ coefficients $(\alpha_k)_{k \ge 0}$ in $\Z_4$
satisfying $\alpha_0 = 1$, $\alpha_1 = 2$ and $\alpha_2 \in \{0,2\}$.
But one can say this:  if the series $f(t)$ possesses an expansion
\reff{def.Stype} with coefficients $(\alpha_k)_{k \ge 0}$ in $R$
{\em and none of these coefficients is a divisor of zero}\/,
then this expansion is unique modulo equivalence
and the algorithm will find it.

The generalization from fields to commutative rings is important
in applications to enumerative combinatorics
\cite{Flajolet_80,Viennot_83,Sokal-Zeng_masterpoly,Sokal_totalpos},
where $R$ is often the ring $\Z[{\bf x}]$ of polynomials
with integer coefficients in some indeterminates ${\bf x} = (x_i)_{i \in I}$.
In particular, the Euler--Gauss recurrence method
applies in an arbitrary commutative ring (with identity element)
and is a useful method for proving continued fractions in this context.

\section{Example 1: From factorials to $\bm{\tfz}$}   \label{sec.examples.1}

Let us now examine Euler's \cite[section~21]{Euler_1760}
derivation of the identity \reff{eq.nfact.contfrac},
which expresses the formal power series $f(t) = \sum_{n=0}^\infty n! \, t^n$
as a regular C-fraction
[that is, \reff{def.Stype} with $p_1 = p_2 = \ldots = 1$]
with coefficients
\be
   \alpha_{2j-1} = j ,\quad \alpha_{2j} = j   \;.
 \label{eq.euler.alpha}
\ee
Euler starts by writing out $f$ through order $t^7$;
he then computes $\alpha_k$ and $f_k$ for $1 \le k \le 7$,
writing each $f_k$ as an explicit ratio $g_k/g_{k-1}$.
It is thus evident that Euler is using what we have here called the
``refined algorithm'' (with $g_{-1} = 1$).
Moreover, Euler writes out each series $g_k$ through order $t^{7-k}$,
to which he appends ``+ etc.'';
clearly he is using the ``finite-$N$ algorithm'' explained in
the preceding section, with $N=7$.
After these computations he says:
\begin{quote}
\small
   And therefore it will become clear, that it will analogously be
$$
   I \,=\, {4x \over 1+K} ,\quad
   K \,=\, {5x \over 1+L} ,\quad
   L \,=\, {5x \over 1+M} ,\quad \hbox{etc.\ to infinity,}
$$
so that the structure of these formulas is easily perceived.
\end{quote}
And he concludes by writing out the continued fraction \reff{eq.nfact.contfrac}
through $\alpha_{13} = 7$ (!), making clear that the intended coefficients
are indeed $\alpha_{2j-1} = j$ and $\alpha_{2j} = j$.

Euler does not give a proof of this final formula
or an explicit expression for the series $g_k$,
but this is not difficult to do.
One approach (the first one I~took)
is to follow Euler and compute the first few coefficients
of the first few $g_k$;
having done this, one can try, by inspecting this finite array of numbers,
to {\em guess}\/ the general formula;
once this has been done, it is not difficult
to {\em prove}\/ the recurrence \reff{eq.gk.recurrence.general}.

But in this case a better approach is available:
namely, compute the {\em full}\/ infinite series $g_k(t)$ for small $k$,
before trying to guess the general formula.
Thus, we begin by setting $g_{-1} = 1$
and $g_0(t) = \sum_{n=0}^\infty n! \, t^n$.
We then use the recurrence \reff{eq.gk.recurrence.general}
[with all $p_i = 1$]
to successively compute $g_1(t)$, $g_2(t)$, \ldots,
extracting at each stage the factor $\alpha_{k+1} t$
that makes $g_{k+1}(t)$ have constant term 1.
After a few steps of this computation,
we may be able to guess the general formulae
for $\alpha_k$ and $g_k(t)$
and then prove the recurrence \reff{eq.gk.recurrence.general}.
Here are the details for this example:

The first step is
\be
   g_0 - g_{-1}
   \;=\;
   \sum_{n=1}^\infty n! \, t^n
   \;=\;
   t \sum_{n=0}^\infty (n+1)! \, t^n
   \;,
\ee
from which we deduce that
$\alpha_1 = 1$ and $g_1(t) = \sum_{n=0}^\infty (n+1)! \, t^n$.
The second step is
\be
   g_1 - g_0
   \;=\;
   \sum_{n=1}^\infty n \, n! \, t^n
   \;=\;
   t \sum_{n=0}^\infty (n+1) \, (n+1)! \, t^n
   \;,
\ee
so that $\alpha_2 = 1$ and $g_2(t) = \sum_{n=0}^\infty (n+1) \, (n+1)! \, t^n$.
Next
\be
   g_2 - g_1
   \;=\;
   \sum_{n=1}^\infty n \, (n+1)! \, t^n
   \;=\;
   2t \sum_{n=0}^\infty (n+1) \, {(n+2)! \over 2} \, t^n
   \;,
\ee
so that $\alpha_3 = 2$ and
$g_3(t) = \sum_{n=0}^\infty (n+1) \, {(n+2)! \over 2} \, t^n$.
And then
\be
   g_3 - g_2
   \;=\;
   \sum_{n=1}^\infty {n(n+1) \over 2} \, (n+1)! \, t^n
   \;=\;
   2t \sum_{n=0}^\infty {(n+1)(n+2) \over 2} \: {(n+2)! \over 2} \, t^n
   \;,
\ee
so that $\alpha_4 = 2$ and
$g_4(t) = \sum_{n=0}^\infty {(n+1)(n+2) \over 2} \: {(n+2)! \over 2} \, t^n$.
At this stage it is not difficult to guess the general formulae
for $\alpha_k$ and $g_k(t)$:
we have $\alpha_{2j-1} = \alpha_{2j} = j$ and
\begin{subeqnarray}
   g_{2j-1}(t) & = & \sum_{n=0}^\infty  \binom{n+j}{n} \binom{n+j-1}{n} \, n! \: t^n
       \\[3mm]
   g_{2j}(t)   & = & \sum_{n=0}^\infty  \binom{n+j}{n}^{\! 2} \, n! \: t^n
 \label{eq.euler.gk}
\end{subeqnarray}
for $j \ge 0$ (as Euler himself may well have known).
Having written down these expressions,
it is then straightforward to verify that they satisfy the recurrence
\be
   g_k(t) - g_{k-1}(t)  \;=\;  \alpha_{k+1} t \, g_{k+1}(t)
   \qquad\hbox{for $k \ge 0$}
 \label{eq.recurrence.gk}
\ee
with the given coefficients $\balpha = (\alpha_k)_{k \ge 1}$.
This completes the proof of \reff{eq.nfact.contfrac}.

In section 26 of the same paper \cite{Euler_1760}, Euler says that
the same method can be applied to the more general series
\reff{eq.xupperfact.contfrac},
which reduces to \reff{eq.nfact.contfrac} when $a=1$;
but he does not provide the details, and he instead
proves \reff{eq.xupperfact.contfrac} by an alternative method.
%
%
Three decades later, however, Euler \cite{Euler_1788}
returned to his original method and presented the details
of the derivation of \reff{eq.xupperfact.contfrac}.\footnote{
   The paper \cite{Euler_1788},
   which is E616 in Enestr\"om's \cite{Enestrom_13} catalogue,
   was apparently presented to the St.~Petersburg Academy in 1776,
   and published posthumously in 1788.
}
By a method similar to the one just shown, one can be led to guess
\be
   \alpha_{2j-1} = a+j-1 ,\quad \alpha_{2j} = j
 \label{eq.euler.alpha.BIS}
\ee
and
\begin{subeqnarray}
   g_{2j-1}(t) & = & \sum_{n=0}^\infty  (a+j)^{\overline{n}} \: \binom{n+j-1}{n} \: t^n
       \\[3mm]
   g_{2j}(t)   & = & \sum_{n=0}^\infty  (a+j)^{\overline{n}} \: \binom{n+j}{n} \: t^n
 \slabel{eq.euler.gk.BIS.2j}
 \label{eq.euler.gk.BIS}
\end{subeqnarray}
where we have used the notation \cite{Knuth_92,Graham_94}
$x^{\overline{n}} = x(x+1) \cdots (x+n-1)$.
The recurrence \reff{eq.recurrence.gk} can once again be easily checked.

We can, in fact, carry this process one step farther,
by introducing an additional parameter $b$.
Let
\be
   \alpha_{2j-1} = a+j-1 ,\quad \alpha_{2j} = b+j-1
 \label{eq.euler.alpha.BISBIS}
\ee
and
\begin{subeqnarray}
   g_{2j-1}(t) & = & \sum_{n=0}^\infty
       {(a+j)^{\overline{n}} \: (b+j-1)^{\overline{n}}  \over n!}  \; t^n
       \\[3mm]
   g_{2j}(t)   & = &
       \sum_{n=0}^\infty
       {(a+j)^{\overline{n}} \: (b+j)^{\overline{n}}  \over n!}  \; t^n
 \label{eq.euler.gk.BISBIS}
\end{subeqnarray}
The recurrence \reff{eq.recurrence.gk} can again be easily checked;
in fact, the reasoning is somewhat more transparent in this greater generality.
We no longer have $g_{-1} = 1$ (unless $b=1$), but no matter;
we can still conclude that $g_0(t)/g_{-1}(t)$
is given by the continued fraction
with coefficients \reff{eq.euler.alpha.BISBIS}.

The series appearing in \reff{eq.euler.gk.BISBIS}
are nothing other than the hypergeometric series $\tfz$, defined by
\be
   \FHYPERbottomzero{2}{a,b}{t}
   \;=\;
   \sum_{n=0}^\infty {a^{\overline{n}} \, b^{\overline{n}}  \over n!}  \; t^n
   \;;
 \label{def.2F0}
\ee
and the recurrence \reff{eq.recurrence.gk} is simply the contiguous relation
\be
   \FHYPERbottomzero{2}{a,b}{t}
   \:-\:
   \FHYPERbottomzero{2}{a,b-1}{t}
   \;=\;
   at \;
   \FHYPERbottomzero{2}{a+1,b}{t}
  \;,
\ee
applied with interchanges $a \leftrightarrow b$ at alternate levels.
We have thus proven the continued fraction for the ratio
of two contiguous hypergeometric series $\tfz$ \cite[section~92]{Wall_48}:
\be
   {\FHYPERbottomzero{2}{a,b}{t}
    \over
    \FHYPERbottomzero{2}{a,b-1}{t}
   }
   \;\:=\;\:
   \cfrac{1}{1 - \cfrac{at}{1 - \cfrac{bt}{1 - \cfrac{(a+1)t}{1- \cfrac{(b+1)t}{1 - \cfrac{(a+2)t}{1 - \cfrac{(b+2)t}{1-\cdots}}}}}}}
   \;\,.
 \label{eq.euler.contfrac.BISBIS.2F0}
\ee
   
At this point let me digress by making three remarks:

\medskip

1) The hypergeometric series \reff{def.2F0}
is of course divergent for all $t \ne 0$
(unless $a$ or $b$ is zero or a negative integer,
 in which case the series terminates).
We can nevertheless give the formula \reff{eq.euler.contfrac.BISBIS.2F0}
an analytic meaning by defining
\be
   F_{a,b}(t)
   \;=\;  
   {1 \over \Gamma(a)} \int\limits_0^\infty
        {e^{-x} \, x^{a-1} \over (1-tx)^b} \: dx
   \;,
 \label{def.Fab}
\ee
which is manifestly an analytic function jointly in $a,b,t$
for $\real a > 0$, $b \in \C$ and $t \in \C \setminus [0,\infty)$;
moreover, its asymptotic expansion at $t=0$
(valid in a sector staying away from the positive real axis)
is the hypergeometric series \reff{def.2F0}.
It~can also be shown that $F_{a,b}(t) = F_{b,a}(t)$
where both sides are defined \cite[p.~277]{Khrushchev_08}.
Furthermore, by integration by parts the definition \reff{def.Fab}
can be extended to arbitrary $a \in \C$.\footnote{
   This is a special case of the more general result that the
   {\em tempered distribution}\/ $x_+^{a-1}/\Gamma(a)$,
   defined initially for $\real a > 0$,
   can be analytically continued to an entire tempered-distribution-valued
   function of $a$ \cite[section~I.3]{Gelfand_64}.
   And this is, in turn, a special case of a spectacular result,
   due to Bernstein and S.I.~Gel'fand \cite{Bernstein_69,Bernstein_72}
   and Atiyah \cite{Atiyah_70},
   on the analytic continuation of distributions of the form
   $P(x_1,\ldots,x_n)^\lambda$ where $P$ is a real polynomial.
   (Here I~digress too far, I~know \ldots\ 
    but this is really beautiful mathematics,
    on the borderline between analysis, algebraic geometry, and algebra:
    see e.g.\ \cite{Bjork_79,Coutinho_95}.)
}
It~can then be shown \cite{Wall_45b}
that the continued fraction on the right-hand side of
\reff{eq.euler.contfrac.BISBIS.2F0}
converges throughout $\C \setminus [0,\infty)$
except possibly at certain isolated points
(uniformly over bounded regions staying away from the isolated points)
and defines an analytic function having these isolated points as poles;
and this analytic function equals $F_{a,b}(t)/F_{a,b-1}(t)$.
(I~know I~had promised to stay away from analysis;
 but this was too beautiful to resist.)

\medskip

2) If we expand the ratio \reff{eq.euler.contfrac.BISBIS.2F0}
as a power series,
\be
   {\FHYPERbottomzero{2}{a,b}{t}
    \over
    \FHYPERbottomzero{2}{a,b-1}{t}
   }
   \;=\;
   \sum_{n=0}^\infty P_n(a,b) \: t^n
   \;,
\ee
it follows easily from the continued fraction that $P_n(a,b)$
is a polynomial of total degree $n$ in $a$ and $b$,
with nonnegative integer coefficients.
It is therefore natural to ask:
What do these nonnegative integers count?

Euler's continued fraction \reff{eq.nfact.contfrac} tells us that
$P_n(1,1) = n!$; and there are $n!$ permutations of an $n$-element set.
It is therefore reasonable to guess that $P_n(a,b)$ enumerates
permutations of an $n$-element set according to some natural
bivariate statistic.
This is indeed the case;
and Dumont and Kreweras \cite{Dumont_88} have identified the statistic.
Given a permutation $\sigma$ of $\{1,2,\ldots,n\}$,
let us say that an index $i \in \{1,2,\ldots,n\}$ is a
\begin{itemize}
   \item \textbfit{record} (or {\em left-to-right maximum}\/)
         if $\sigma(j) < \sigma(i)$ for all $j < i$
      [note in particular that the index 1 is always a record];
   \item \textbfit{antirecord} (or {\em right-to-left minimum}\/)
         if $\sigma(j) > \sigma(i)$ for all $j > i$
      [note in particular that the index $n$ is always an antirecord];
   \item \textbfit{exclusive record} if it is a record and not also
         an antirecord;
   \item \textbfit{exclusive antirecord} if it is an antirecord and not also
      a record.
\end{itemize}
Dumont and Kreweras \cite{Dumont_88} then showed that
\be
   P_n(a,b)  \;=\;  \sum_{\sigma \in \Sym_n} a^{\rec(\sigma)} b^{\earec(\sigma)}
 \label{eq.Pnab}
\ee
where $\rec(\sigma)$ [resp.\ $\earec(\sigma)$]
is the number of records (resp.\ exclusive antirecords) in~$\sigma$.
Some far-reaching generalizations of this result can be found in
\cite{Sokal-Zeng_masterpoly}.

\medskip

3)  Euler also observed \cite[section~29]{Euler_1760} that the case $a = \half$
of \reff{eq.xupperfact.contfrac} leads to
\be
   \sum_{n=0}^\infty (2n-1)!! \: t^n
   \;=\;
   \cfrac{1}{1 - \cfrac{1t}{1 - \cfrac{2t}{1 - \cfrac{3t}{1- \cfrac{4t}{1- \cdots}}}}}
 \label{eq.2n-1semifact.contfrac}
\ee
with coefficients $\alpha_k = k$.
Since $(2n-1)!! = (2n)!/(2^n n!)$
is the number of perfect matchings of a $2n$-element set
(i.e.\ partitions of the $2n$ objects into $n$ pairs),
it is natural to seek generalizations of \reff{eq.2n-1semifact.contfrac}
that enumerate perfect matchings according to some combinatorially
interesting statistics.
Some formulae of this type can be found in
\cite{Dumont_86,Sokal-Zeng_masterpoly}.
The proofs use the bijective method to be discussed
at the end of Section~\ref{sec.combinatorial};
I~don't know whether results of this complexity
can be proven by the Euler--Gauss recurrence method.

\bigskip

This is by no means the end of the matter:
by an argument similar to the one we have used for $\tfz$,
Gauss \cite{Gauss_1813} found in 1812
a continued fraction for the ratio of
two contiguous hypergeometric functions $\tfo$.
Moreover, the formula for $\tfz$,
as well as analogous formulae for ratios of $\ofo$, $\ofz$ or $\zfo$,
can be deduced from Gauss' formula by specialization or taking limits.
In fact, one of the special cases of the $\zfo$ formula
is Lambert's continued fraction \reff{eq.lambert} for the tangent function.
See \cite[Chapter~XVIII]{Wall_48} for details.\footnote{
   Laczkovich \cite{Laczkovich_97} and Wallisser \cite{Wallisser_00}
   give nice elementary proofs of the continued fraction for $\zfo$,
   using the Euler--Gauss recurrence method.
   As Wallisser \cite[p.~525]{Wallisser_00} points out,
   this argument is due to Legendre
   \cite[Note~IV, pp.~320--322]{Legendre_1800}.
}
There is also a nice explanation at \cite{wikipedia},
which makes clear the general principle
of the Euler--Gauss recurrence method:
any recurrence of the form \reff{eq.gk.recurrence.general}
for a sequence $(g_k)_{k \ge -1}$ of series with constant term 1
leads to a continued-fraction representation \reff{def.fk}
for the ratios $f_k = g_k/g_{k-1}$.
In Sections~\ref{sec.Jtype} and \ref{sec.cfgen}
we will see even more general versions of this principle.

\section{Example 2: Bell polynomials}   \label{sec.examples.2}

Here is an example from enumerative combinatorics.
The \textbfit{Bell number} $B_n$ is, by~definition,
the number of partitions of an $n$-element set into nonempty blocks;
by~convention we set $B_0 = 1$.
The \textbfit{Stirling subset number}
(also called Stirling number of the second kind)
$\stirlingsubset{n}{k}$ is, by~definition,
the number of partitions of an $n$-element set into $k$ nonempty blocks;
for $n=0$ we make the convention $\stirlingsubset{0}{k} = \delta_{k0}$.
The Stirling subset numbers satisfy the recurrence
\be
   \stirlingsubset{n}{k}
   \;=\;
   k \, \stirlingsubset{n-1}{k}
   \:+\:
   \stirlingsubset{n-1}{k-1}
   \qquad\hbox{for $n \ge 1$}
 \label{eq.stirling.recurrence}
\ee
with initial conditions $\stirlingsubset{0}{k} = \delta_{k0}$
and $\stirlingsubset{n}{-1} = 0$.
[{\sc Proof:}
 Consider a partition $\pi$ of the set $[n] \eqdef \{1,\ldots,n\}$
 into $k$ nonempty blocks, and ask where the element $n$ goes.
 If the restriction of $\pi$ to $[n-1]$ has $k$ blocks,
 then $n$ can be adjoined to any one of those $k$ blocks.
 If the restriction of $\pi$ to $[n-1]$ has $k-1$ blocks,
 then $n$ must be a singleton in $\pi$.
 These two cases give the two terms on the right-hand side of
 \reff{eq.stirling.recurrence}.]

Now define the \textbfit{Bell polynomials}
\be
   B_n(x)  \;=\;  \sum_{k=0}^n \stirlingsubset{n}{k} \, x^k
\ee
and their homogenized version
\be
   B_n(x,y)  \;=\;  y^n B_n(x/y)
             \;=\;  \sum_{k=0}^n \stirlingsubset{n}{k} \, x^k y^{n-k}
   \;,
\ee
so that $B_n = B_n(1) = B_n(1,1)$.
Then the ordinary generating function
\be
   \scrb_{x,y}(t)  \;=\;  \sum_{n=0}^\infty B_n(x,y) \, t^n
\ee
turns out to have a beautiful continued fraction:
\be
   \scrb_{x,y}(t)
   \;=\;
   \cfrac{1}{1 - \cfrac{xt}{1 - \cfrac{yt}{1 - \cfrac{xt}{1- \cfrac{2yt}{1 - \cfrac{xt}{1 - \cfrac{3yt}{1-\cdots}}}}}}}
   \label{def.contfrac.Bxy}
\ee
with coefficients $\alpha_{2k-1} = x$ and $\alpha_{2k} = ky$.

Once again we can guess the continued fraction, and then prove it,
by the Euler--Gauss recurrence method.
Take $g_{-1} = 1$ and $g_0(t) = \scrb_{x,y}(t)$,
and use the recurrence \reff{eq.recurrence.gk}
to successively compute $g_1(t)$, $g_2(t)$, \ldots,
extracting at each stage the factor $\alpha_{k+1} t$
that makes $g_{k+1}(t)$ have constant term 1.
This computation is left as an exercise for the reader;
by the stage $g_6$ (if not earlier)
the reader should be able to guess the general formulae
for $g_{2j-1}(t)$ and $g_{2j}(t)$.
(In order not to spoil the fun, the answer is given in the Appendix.)
Once one has the formulae for $g_k(t)$,
it is then easy to verify the recurrence \reff{eq.recurrence.gk}
with the given coefficients $\balpha$
by using the recurrence \reff{eq.stirling.recurrence}
for the Stirling subset numbers
together with the Pascal recurrence for the binomial coefficients.

\bigskip

{\bf Remarks.}
I am not sure who first derived the continued fraction
\reff{def.contfrac.Bxy} for the Bell polynomials,
or its specialization to $x=y=1$ for the Bell numbers.  
An associated continued fraction\footnote{
   In the terminology of combinatorialists, a J-fraction.
}
that is equivalent by contraction
\cite[p.~21]{Wall_48} \cite[p.~V-31]{Viennot_83}
to \reff{def.contfrac.Bxy}
was found for the case $x=y=1$
by Touchard \cite[section~4]{Touchard_56} in 1956,
and for the general case
by Flajolet \cite[Theorem~2(ia)]{Flajolet_80} in 1980.
Flajolet's proof was combinatorial,
using ideas that will be explained in Section~\ref{sec.combinatorial}.
Flajolet also observed \cite[pp.~141--142]{Flajolet_80}
that this associated continued fraction is implicit
in the three-term recurrence relation for the Poisson--Charlier polynomials
\cite[p.~25, Exercise~4.10]{Chihara_78};
see \cite{Chihara_78,Viennot_83,Zeng_21} for the general connection between
continued fractions and orthogonal polynomials.
The continued fraction \reff{def.contfrac.Bxy}
can also be derived directly from a functional equation satisfied by
$\scrb_{x,y}(t)$:
this elegant method is due to the late Dominique Dumont \cite{Dumont_89};
see also \cite[proof of Lemma~3]{Zeng_95} for some $q$-generalizations.
I have not seen the elementary derivation by the Euler--Gauss recurrence method
anywhere in the literature, but it is probably not new.

See \cite{Kasraoui_06,Sokal-Zeng_masterpoly} for some generalizations
of this continued fraction, which enumerate set partitions
with respect to a larger set of simultaneous statistics;
these formulae are proven by the bijective method to be discussed
at the end of Section~\ref{sec.combinatorial}.
\myendremark

\section{Example 3: \hbox{Some $\bm{q}$-continued fractions of Ramanujan}}
    \label{sec.examples.3}

Next I would like to show, following Bhatnagar \cite{Bhatnagar_14},
how the Euler--Gauss recurrence method can be used to give simple proofs
of some continued fractions of Ramanujan.
We use the standard notation for $q$-shifted factorials,
\be
   (a;q)_n  \;=\;  \prod_{j=0}^{n-1}  (1- aq^j)
\ee
for integers $n \ge 0$;
here $a$ and $q$ are to be interpreted as algebraic indeterminates.

\bigskip

{\bf The Rogers--Ramanujan continued fraction.}
Rogers \cite[p.~328, eq.~(4)]{Rogers_1894} proved in 1894
the following beautiful continued fraction,
which was later rediscovered and generalized by Ramanujan
\cite{Watson_29} \cite[p.~30, Entry~15 and Corollary]{Berndt_91}:
\be
   {\displaystyle \;\sum_{n=0}^\infty  {q^{n^2} \over (q;q)_n} \: t^n \;
    \over
    \displaystyle \;\sum_{n=0}^\infty  {q^{n(n-1)} \over (q;q)_n} \: t^n \;
   }
   \;\:=\;\:
   \cfrac{1}{1 + \cfrac{t}{1 + \cfrac{qt}{1 + \cfrac{q^2 t}{1+ \cfrac{q^3 t}{1 + \cdots}}}}}
 \label{eq.cfrac.RR}
\ee
with coefficients $\alpha_k = - q^{k-1}$.
The proof by the Euler--Gauss recurrence method is extraordinarily easy.
Define
\be
   g_k(t)  \;=\;  \sum_{n=0}^\infty  {q^{n(n+k)} \over (q;q)_n} \: t^n
   \qquad\hbox{for $k \ge -1$}
   \;,
 \label{def.gk.RR}
\ee
so that the left-hand side of \reff{eq.cfrac.RR} is indeed $g_0/g_{-1}$.
Then compute
\begin{subeqnarray}
  g_k - g_{k-1}
  & = &
  \sum_{n=0}^\infty  {q^{n(n+k-1)} \, (q^n - 1) \over (q;q)_n} \: t^n
     \\[2mm]
  & = &
  - \sum_{n=1}^\infty  {q^{n(n+k-1)} \over (q;q)_{n-1}} \: t^n
     \\[2mm]
  & = &
  - \sum_{n=0}^\infty  {q^{(n+1)(n+k)} \over (q;q)_n} \: t^{n+1}
     \\[2mm]
  & = &
  - q^k t \sum_{n=0}^\infty  {q^{n(n+k+1)} \over (q;q)_n} \: t^n
     \\[2mm]
  & = &
  \alpha_{k+1} t \, g_{k+1}
  \;,
 \label{eq.RR.computation}
\end{subeqnarray}
which completes the proof
(see also \cite[eqns.~(4.43)/(4.44)]{Askey_89} \cite{Bhatnagar_15}).

In terms of the \textbfit{Rogers--Ramanujan function}
\be
   R(t,q)  \;=\;  \sum_{n=0}^\infty  {q^{n(n-1)} \over (q;q)_n} \: t^n
   \;,
\ee
we have $g_k(t) = R(q^{k+1}t,q)$;
the left-hand side of \reff{eq.cfrac.RR} is
$f_0(t) = R(qt,q) / R(t,q)$,
and more generally we have $f_k(t) = R(q^{k+1}t,q) / R(q^k t,q)$.
It is worth remarking that the Rogers--Ramanujan function
arises in a two-variable identity
due to Ramanujan and Rogers \cite{Ramanujan_19}
from which the famous one-variable Rogers--Ramanujan identities
\cite[Chapter~7]{Andrews_76} \cite{Andrews_89,Sills_18} can be deduced.
The Rogers--Ramanujan function has also been studied
as an entire function of $t$ for $|q| < 1$ \cite{Andrews_05a}.

In fact, Ramanujan \cite[p.~30, Entry~15]{Berndt_91}
gave a generalization of \reff{eq.cfrac.RR} with an additional free parameter;
this result can be rewritten \cite[p.~57, Exercise]{Bhatnagar_14} as
\be
   {\displaystyle \;\sum_{n=0}^\infty  {q^{n^2} \over (q;q)_n \, (a;q)_n}
                                                                      \: t^n \;
    \over
    \displaystyle \;\sum_{n=0}^\infty  {q^{n(n-1)} \over (q;q)_n \, (a;q)_n}
                                                                      \: t^n \;
   }
   \;\:=\;\:
   \cfrac{1}{1 + \cfrac{\disfrac{1}{1-a}\,t}{1 + \cfrac{\disfrac{q}{(1-a)(1-aq)}\,t}{1 + \cfrac{\disfrac{q^2}{(1-aq)(1-aq^2)}\,t}{1+ \cfrac{\disfrac{q^3}{(1-aq^2)(1-aq^3)}\,t}{1 + \cdots}}}}}
 \label{eq.cfrac.RRbis}
\ee
with coefficients
\be
   \alpha_1 \;=\;  - \, {1 \over 1-a} \:,\qquad
   \alpha_k \;=\;  - \, {q^{k-1} \over (1-aq^{k-2})(1-aq^{k-1})}
       \quad\hbox{for $k \ge 2$}
   \;.
\ee
(Note the difference in form between $\alpha_1$ and the remaining coefficients:
 one factor in the denominator versus two.)
This result can be derived by a slight generalization of
the computation \reff{eq.RR.computation}, using
\begin{subeqnarray}
   g_{-1}(t)
   & = &
   \sum_{n=0}^\infty  {q^{n(n-1)} \over (q;q)_n \, (a;q)_n} \: t^n
          \\[3mm]
   g_k(t)
   & = &
   \sum_{n=0}^\infty  {q^{n(n+k)} \over (q;q)_n \, (aq^k;q)_n} \: t^n
      \qquad\hbox{for $k \ge 0$}
 \label{def.gk.RRbis}
\end{subeqnarray}
(Note the corresponding difference between $k=-1$ and $k \ge 0$.)
The proof, which is not difficult, is left as an exercise for the reader.

On the other hand, there is a variant of \reff{eq.cfrac.RRbis}
that is even simpler.  Namely, use $(aq;q)_n$ instead of $(a;q)_n$
in the numerator of the left-hand side (but not the denominator);
then we have
\be
   {\displaystyle \;\sum_{n=0}^\infty  {q^{n^2} \over (q;q)_n \, (aq;q)_n}
                                                                      \: t^n \;
    \over
    \displaystyle \;\sum_{n=0}^\infty  {q^{n(n-1)} \over (q;q)_n \, (a;q)_n}
                                                                      \: t^n \;
   }
   \;\:=\;\:
   \cfrac{1}{1 + \cfrac{\disfrac{1}{(1-a)(1-aq)}\,t}{1 + \cfrac{\disfrac{q}{(1-aq)(1-aq^2)}\,t}{1 + \cfrac{\disfrac{q^2}{(1-aq^2)(1-aq^3)}\,t}{1+ \cfrac{\disfrac{q^3}{(1-aq^3)(1-aq^4)}\,t}{1 + \cdots}}}}}
 \label{eq.cfrac.RRbis2}
\ee
with coefficients
\be
   \alpha_k \;=\;  - \, {q^{k-1} \over (1-aq^{k-1})(1-aq^{k})}
   \;.
\ee
Now there is no difference between the first step and the rest,
and we can use the single formula
\be
   g_k(t)
   \;=\;
   \sum_{n=0}^\infty  {q^{n(n+k)} \over (q;q)_n \, (aq^{k+1};q)_n} \: t^n
 \label{def.gk.RRbis2}
\ee
for all $k \ge -1$.
In terms of the basic hypergeometric series $\phiHyper{r}{s}$
defined by \cite[p.~4]{Gasper_04}
\be
   \phiHYPER{r}{s\!}{a_1,\ldots,a_r}{b_1,\ldots,b_s}{q}{t}
   \;=\;
   \sum_{n=0}^\infty
   {(a_1;q)_n \, (a_2;q)_n \,\cdots\, (a_r;q)_n
    \over
    (b_1;q)_n \, (b_2;q)_n \,\cdots\, (b_s;q)_n \, (q;q)_n
   }
   \:
   \Bigl(\! (-1)^n q^{n(n-1)/2} \!\Bigr) ^{\! s+1-r}
   \:
   t^n
   \;,
 \label{def.rphis}
\ee
the left-hand side of \reff{eq.cfrac.RRbis2} is
$
   \displaystyle \phiHYPER{0}{1\!}{\hbox{---}}{aq}{q}{qt}
   \biggl/\!
    \displaystyle \phiHYPER{0}{1\!}{\hbox{---}}{a}{q}{t}
$\,,
and the continued fraction \reff{eq.cfrac.RRbis2}
can alternatively be derived as a limiting case of Heine's
\cite{Heine_1847} \cite[p.~395]{Cuyt_08}
continued fraction for ratios of contiguous $\phiHyper{2}{1}$.




\bigskip

{\bf The partial theta function.}
The function
\be
   \Theta_0(t,q)  \;=\; \sum_{n=0}^\infty q^{n(n-1)/2} \, t^n
\ee
is called the \textbfit{partial theta function}
\cite[Chapter~13]{Andrews-Berndt_05} \cite[Chapter~6]{Andrews-Berndt_09}
\cite{Andrews_05b,Sokal_partialtheta}
because of its resemblance to the ordinary theta function,
in which the sum runs down to ${n=-\infty}$.
A continued-fraction expansion for the partial theta function
was discovered by Eisenstein \cite{Eisenstein_1844a,Eisenstein_1844b} in 1844
and rediscovered by Ramanujan \cite[pp.~27--29, Entry~13]{Berndt_91}
(see also \cite{Ramanathan_87,Folsom_06}).
It reads
\be
   \sum_{n=0}^\infty q^{n(n-1)/2} \, t^n
   \;=\;
   \cfrac{1}{1 - \cfrac{t}{1 - \cfrac{(q-1)t}{1 - \cfrac{q^2 t}{1- \cfrac{q(q^2-1)t}{1- \cfrac{q^4 t}{1- \cfrac{q^2 (q^3-1)t}{1- \cdots}}}}}}}
 \label{eq.partialtheta.contfrac}
\ee
with coefficients
\be
   \alpha_{2j-1} \;=\;  q^{2j-2} ,\qquad
   \alpha_{2j}   \;=\;  q^{j-1} (q^j - 1)
   \;.
 \label{eq.partialtheta.alphas}
\ee

Once again we can guess the continued fraction, and then prove it,
by the Euler--Gauss recurrence method with $g_{-1} = 1$;
but here it is a bit trickier than in the previous examples
to guess the coefficients $\balpha$ and the series $g_k(t)$.
The computation is once again left as an exercise for the reader;
by the stage $g_6$ it should become clear that
the coefficients $\balpha$ are given by \reff{eq.partialtheta.alphas}
and the series $g_k(t)$ by
\begin{subeqnarray}
   g_{2j-1}(t)
   & = &
   \sum_{n=0}^\infty \qbinom{n+j-1}{n}{\!q} \; q^{n(n+2j-1)/2} \: t^n
       \\[3mm]
   g_{2j}(t)
   & = &
   \sum_{n=0}^\infty \qbinom{n+j}{n}{\!q} \; q^{n(n+2j-1)/2} \: t^n
 \label{eq.partialtheta.gk}
\end{subeqnarray}
where $\qbinom{n}{k}{q}$ denotes the \textbfit{$\bm{q}$-binomial coefficient}
\be
   \qbinom{n}{k}{\!q}
   \;=\;
   {(q;q)_n  \over  (q;q)_k \, (q;q)_{n-k}}
   \;.
 \label{def.qbinom.bis}
\ee
The right-hand side of \reff{def.qbinom.bis}
looks like a rational function of $q$,
but it is a nontrivial fact (though not terribly difficult to prove)
that $\qbinom{n}{k}{q}$ is in fact a {\em polynomial}\/ in $q$,
with nonnegative integer coefficients that have a nice combinatorial
interpretation \cite[Theorem~3.1]{Andrews_76}.
The $q$-binomial coefficients satisfy two ``dual'' $q$-generalizations
of the Pascal recurrence:
\begin{eqnarray}
   \qbinom{n}{k}{\!q}
   & = &
   \qbinom{n-1}{k}{\!q} \,+\, q^{n-k} \qbinom{n-1}{k-1}{\!q}
       \quad\hbox{for $n \ge 1$}
 \label{eq.binom.recursion1}  \\[3mm]
   \qbinom{n}{k}{\!q}
   & = &
   q^k \qbinom{n-1}{k}{\!q} \,+\, \qbinom{n-1}{k-1}{\!q}
       \qquad\hbox{for $n \ge 1$}
 \label{eq.binom.recursion2}
\end{eqnarray}
(Of course, it follows immediately from either of these recurrences
that $\qbinom{n}{k}{q}$ is a polynomial in $q$,
with nonnegative integer coefficients.)
Using the recurrence \reff{eq.binom.recursion1},
it is now straightforward to verify the
Euler--Gauss recurrence \reff{eq.recurrence.gk}
for the given $\balpha$ and $g_k$.
This completes the proof of \reff{eq.partialtheta.contfrac}.

A different (but also simple) proof of \reff{eq.partialtheta.contfrac}
is given in \cite[pp.~27--28, Entry~13]{Berndt_91}.
A more general continued fraction can be found in Ramanujan's lost notebook:
see \cite[Section~6.2]{Andrews-Berndt_05}.


The reader is referred to Bhatnagar's beautiful survey articles
\cite{Bhatnagar_14,Bhatnagar_22}
for derivations of many other continued fractions of Ramanujan
by the Euler--Gauss recurrence method (among other methods).
See also \cite{Gu_11} for a cornucopia of related results.

\section{Expansion in the form \reff{def.Jtype}}   \label{sec.Jtype}

Let us now consider expansion in the form \reff{def.Jtype},
which generalizes the C-fraction \reff{def.Stype} and reduces to it when
$M_1 = M_2 = \ldots = 0$.
Here we consider the integers $M_i \ge 0$ to be pre-specified,
while the integers $p_i \ge M_i + 1$ are chosen by the algorithm.

Since the treatment closely parallels that of \reff{def.Stype},
I~will be brief and stress only the needed modifications.
It is convenient to use the abbreviation
\be
   \Delta_i(t)
   \;=\;
   \sum\limits_{j=1}^{M_i} \delta_i^{(j)} t^j
\ee
for the ``additive'' coefficient in \reff{def.Jtype};
it is a polynomial of degree $\le M_i$ in $t$, with zero constant term.

As usual we define
\be
   f_k(t)  \;=\;
   \cfrac{1}{1 - \Delta_{k+1}(t) - \cfrac{\alpha_{k+1} t^{p_{k+1}}}{1 - \Delta_{k+2}(t) - \cfrac{\alpha_{k+2} t^{p_{k+2}}}{1 - \cdots}}}
   \label{def.fk.Jtype}
\ee
and observe that $f(t) = \alpha_0 f_0(t)$ and
\be
   f_k(t)
   \;=\;
   {1 \over 1 \,-\, \Delta_{k+1}(t) \,-\, \alpha_{k+1} t^{p_{k+1}} \, f_{k+1}(t)}
   \qquad\hbox{for $k \ge 0$}  \;.
 \label{eq.recurrence.fk.Jtype}
\ee
The primitive algorithm is then:
%
%
\begin{framed}[0.82\textwidth]
\noindent
{\bf Primitive algorithm.}

\medskip
\noindent
1.  Set $\alpha_0 = a_0 = [t^0] \, f(t)$ and
      $f_0(t) = \alpha_0^{-1} f(t)$.

\medskip
\noindent
2.  For $k=1,2,3,\ldots$, do:
\begin{itemize}
   \item[(a)]  Set $\Delta_k(t)$ equal to the expansion of
       $1 - f_{k-1}(t)^{-1}$ \hbox{through order $t^{M_k}$.}
   \item[(b)]  If $1 - f_{k-1}(t)^{-1} = \Delta_k(t)$,
       set $\alpha_k = 0$ and terminate.
   \item[(c)]  Otherwise, let $p_k$ be the smallest index $n > M_k$
       such that $[t^n] \, f_{k-1}(t)^{-1} \neq 0$;
       set $\alpha_k = - [t^{p_k}] \, f_{k-1}(t)^{-1}$;
       and set
\be
   f_k(t)  \;=\;  \alpha_k^{-1} t^{-p_k} \biggl( 1 \,-\, {1 \over f_{k-1}(t)}
                                                   \,-\, \Delta_k(t)
                                         \biggr)
   \;.
 \label{eq.alg.fk.Jtype}
\ee
\end{itemize}
\vspace*{-5mm}
\end{framed}

\bigskip

{\bf Historical remark.}
The case $M_1 = M_2 = \ldots =1$ of the primitive algorithm
was proposed in 1772 by Lagrange \cite{Lagrange_1772}.
See Brezinski \cite[pp.~119--120]{Brezinski_91}
and especially Galuzzi \cite{Galuzzi_95}
for further discussion of this work.
\myendremark

\bigskip

Let us now discuss the refined algorithm,
passing immediately to the generalized version
in which $g_{-1}$ is an arbitrary series with constant term 1.
The series $(g_k)_{k \ge 0}$ are therefore defined by \reff{def.gk.generalized},
so that $f_k = g_k/g_{k-1}$ as before.
Then the nonlinear recurrence \reff{eq.recurrence.fk.Jtype} for the $(f_k)$
becomes the linear recurrence
\be
   g_k(t) - g_{k-1}(t)
   \;=\;
   \Delta_{k+1}(t) g_k(t) \,+\, \alpha_{k+1} t^{p_{k+1}} g_{k+1}(t)
 \label{eq.recurrence.gk.Jtype}
\ee
for the $(g_k)$.
The occurrence here of the term $\Delta_{k+1} g_k$ means that
division of power series is now required in order to determine $\Delta_{k+1}$;
but this division need only be exact through order $t^{M_{k+1}}$,
which is not onerous if $M_{k+1}$ is small.
Rewriting the algorithm in terms of $(g_k)_{k \ge -1}$, we have:

\begin{framed}[0.88\textwidth]
\noindent
{\bf Refined algorithm.}

\medskip
\noindent
1. Choose any formal power series $g_{-1}(t)$ with constant term~1;
      \linebreak
      then set $\alpha_0 = a_0 = [t^0] \, f(t)$ and
      $g_0(t) = \alpha_0^{-1} g_{-1}(t) f(t)$.

\medskip
\noindent
2. For $k=1,2,3,\ldots$, do:
\begin{itemize}
   \item[(a)]  Set $\Delta_k(t)$ equal to the expansion of
       $1 - g_{k-2}(t)/g_{k-1}(t)$ \hbox{through order $t^{M_k}$.}
   \item[(b)]  If $g_{k-1}(t) - g_{k-2}(t) - \Delta_k(t) g_{k-1}(t) = 0$,
      set $\alpha_k = 0$ and terminate.
   \item[(c)]  Otherwise, let $p_k$ be the smallest index $n$
      (necessarily $> M_k$) such that
      $[t^n] \, \bigl( g_{k-1}(t) - g_{k-2}(t) - \Delta_k(t) g_{k-1}(t) \bigr)
       \neq 0$;
      \hspace*{2cm}\linebreak
      set $\alpha_k = [t^{p_k}] \,
                \bigl( g_{k-1}(t) - g_{k-2}(t) - \Delta_k(t) g_{k-1}(t) \bigr)$;
       and set
\be
   g_k(t)  \;=\;  \alpha_k^{-1} t^{-p_k}
               \bigl( g_{k-1}(t) - g_{k-2}(t) - \Delta_k(t) g_{k-1}(t) \bigr)
   \;.
 \label{eq.alg.gk.Jtype}
\ee
\end{itemize}
\vspace*{-5mm}
\end{framed}

\bigskip

We can also run this algorithm in reverse,
leading to a generalization of the Euler--Gauss recurrence method
as presented in \reff{eq.gk.recurrence.general}.
Suppose that we have a sequence $(g_k)_{k \ge -1}$
of formal power series with constant term 1,
which satisfy a recurrence of the form
\be
   g_{k}(t) - g_{k-1}(t)  \;=\;
    \Delta_{k+1}(t) \, g_{k}(t) \:+\: A_{k+1}(t) \, g_{k+1}(t)
   \qquad\hbox{for $k \ge 0$}
 \label{eq.gk.recurrence.Jtype}
\ee
where the $\Delta_k(t)$ and $A_k(t)$ 
are formal power series with zero constant term.
(We need not assume that $g_{-1} = 1$,
 nor that $\Delta_k(t)$ is a polynomial,
 nor that $A_k(t)$ is simply a monomial $\alpha_k t^{p_k}$.)
Dividing by $g_{k}$ and defining $f_k = g_k/g_{k-1}$,
we have
\be
   f_{k}(t)
   \;=\;
   {1 \over 1 \,-\, \Delta_{k+1}(t) \,-\, A_{k+1}(t) \, f_{k+1}(t)}
   \qquad\hbox{for $k \ge 0$}  \;,
 \label{eq.recurrence.fk.Jtype.BIS}
\ee
which by iteration yields the continued-fraction expansions
\be
   f_k(t)  \;=\;
   \cfrac{1}{1 - \Delta_{k+1}(t) - \cfrac{A_{k+1}(t)}{1 - \Delta_{k+2}(t) - \cfrac{A_{k+2}(t)}{1 - \cdots}}}
   \;.
   \label{def.fk.Jtype.BIS}
\ee
When $\Delta_k(t)$ is a polynomial of degree $\le M_k$
and $A_k(t) = \alpha_k t^{p_k}$,
this reduces to \reff{def.fk.Jtype}.
This method was used by Rogers \cite[p.~76]{Rogers_07} in 1907
to obtain expansions as an associated continued fraction
(i.e.\ $M_1 = M_2 = \ldots =1$ and $p_1 = p_2 = \ldots = 2$)
for the Laplace transforms of the Jacobian elliptic functions sn and cn
(see also \cite[p.~237]{Flajolet_89}).
Some spectacular extensions of these results, using the same method,
were given in the early 2000s by Milne \cite[section~3]{Milne_02}
and Conrad and Flajolet \cite{Conrad_02,Conrad_06}.
On the other hand,
the special case $\Delta_k(t) = \delta_k t$ and $A_k(t) = \alpha_k t$
is also important, and is called a T-fraction
\cite{Thron_48,Roblet_96,Sokal_totalpos,latpath_SRTR,Elvey-Sokal_wardpoly}.

\section{Expansion in the form \reff{def.cfgen}}   \label{sec.cfgen}

The continued-fraction schema \reff{def.cfgen} is so general that the expansion
of a given series $f(t)$ in this form is far from unique.
Indeed, the series $\Delta_k(t)$ can be chosen completely arbitrarily
(with zero constant term),
while the $A_k(t)$ need only have the correct leading terms
and are otherwise also completely arbitrary.
Let us define as usual
\be
   f_k(t)  \;=\;
   \cfrac{1}{1 - \Delta_{k+1}(t) - \cfrac{A_{k+1}(t)}{1 - \Delta_{k+2}(t) - \cfrac{A_{k+2}(t)}{1 - \cdots}}}
   \qquad\hbox{for $k \ge 0$}  \;;
   \label{def.fk.cfgen}
\ee
these are formal power series with constant term 1,
which satisfy $f(t) = A_0(t) \, f_0(t)$ and
\be
   f_k(t)
   \;=\;
   {1 \over 1 \,-\, \Delta_{k+1}(t) \,-\, A_{k+1}(t) \, f_{k+1}(t)}
   \qquad\hbox{for $k \ge 0$}  \;.
 \label{eq.recurrence.fk.cfgen}
\ee
The procedure for finding a continued-fraction expansion
of a given series $f(t)$ in the form \reff{def.cfgen}
--- I~am reluctant to call it an ``algorithm'',
as it now involves so many arbitrary choices ---
is then as follows:

\noindent
\begin{framed}[0.88\textwidth]
\noindent
{\bf Primitive procedure.}

\medskip
\noindent
1.  Let $A_0(t)$ be any formal power series
      having the same leading term as $f(t)$;
      and set $f_0(t) = A_0(t)^{-1} f(t)$.

\medskip
\noindent
2.  For $k=1,2,3,\ldots$, do:
\begin{itemize}
   \item[(a)] Let $\Delta_k(t)$ be any formal power series with zero
     constant term.
   \item[(b)]  If $1 - f_{k-1}(t)^{-1} = \Delta_k(t)$, set $A_k(t) = 0$
      and terminate.
   \item[(c)]  Otherwise, let $p_k$ be the smallest index $n$
       such that
       \linebreak
       $[t^n] \, [1 - f_{k-1}(t)^{-1} - \Delta_k(t)] \neq 0$;
       set $\alpha_k = [t^{p_k}] \, [1 - f_{k-1}(t)^{-1} - \Delta_k(t)]$;
       let $A_k(t)$ be any formal power series
       with leading term $\alpha_k t^{p_k}$;
       and set
\be
   f_k(t)  \;=\;  A_k(t)^{-1} \biggl( 1 \,-\, {1 \over f_{k-1}(t)}
                                        \,-\, \Delta_k(t)
                                         \biggr)
   \;.
 \label{eq.alg.fk.cfgen}
\ee
\end{itemize}
\end{framed}

The corresponding refined procedure is now left as an exercise
for the reader; it is a minor modification of the one presented
in the preceding section.
And the corresponding generalization of the Euler--Gauss recurrence method
was already discussed in that section.



\section{Combinatorial interpretation}   \label{sec.combinatorial}

A combinatorial interpretation of continued fractions
in terms of lattice paths was given
in a seminal 1980 paper by the late Philippe Flajolet \cite{Flajolet_80};
we review it here, and then show how it can be used to interpret
the series $(f_k)_{k \ge 0}$ and $(g_k)_{k \ge 0}$
arising in our algorithm.

A \textbfit{Motzkin path}
is a path in the upper half-plane $\Z \times \N$,
starting and ending on the horizontal axis,
using steps $(1,1)$ [``rise''], $(1,0)$ [``level step'']
and $(1,-1)$ [``fall''].
More generally, a \textbfit{Motzkin path at level~$\bm{k}$}
is a path in $\Z \times \N_{\ge k}$,
starting and ending at height $k$, using these same steps.
We denote by $\scrm_{k \to k}$
the set of all Motzkin paths at level~$k$ that start at $(0,k)$.
We stress that a Motzkin path must always stay on or above the horizontal axis,
and that a Motzkin path at level~$k$ must always stay at height $\ge k$.
A Motzkin path is called a \textbfit{Dyck path} if it has no level steps;
obviously a Dyck path must have even length.

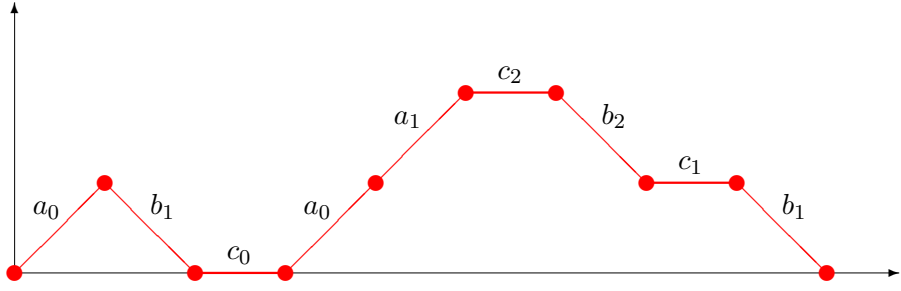
\begin{figure}[t]
\centering
\vspace*{3cm}
\hspace*{-11.5cm}
\begin{picture}(12,5)(0,0)
\setlength{\unitlength}{1.2cm}
\linethickness{0.1mm}
\put(0,0){\vector(1,0){9.8}}
\put(0,0){\vector(0,1){3}}
\linethickness{0.3mm}
\put(0,0){\red{\line(1,1){1}}}\put(0.2,0.65){\small $a_0$}
\put(1,1){\red{\line(1,-1){1}}}\put(1.5,0.65){\small $b_1$}
\put(2,0){\red{\line(1,0){1}}}\put(2.35,0.15){\small $c_0$}
\put(3,0){\red{\line(1,1){1}}}\put(3.2,0.65){\small $a_0$}
\put(4,1){\red{\line(1,1){1}}}\put(4.2,1.65){\small $a_1$}
\put(5,2){\red{\line(1,0){1}}}\put(5.35,2.15){\small $c_2$}
\put(6,2){\red{\line(1,-1){1}}}\put(6.5,1.65){\small $b_2$}
\put(7,1){\red{\line(1,0){1}}}\put(7.35,1.15){\small $c_1$}
\put(8,1){\red{\line(1,-1){1}}}\put(8.5,0.65){\small $b_1$}
\put(0,0){\red{\circle*{.17}}}
\put(1,1){\red{\circle*{.17}}}
\put(2,0){\red{\circle*{.17}}}
\put(3,0){\red{\circle*{.17}}}
\put(4,1){\red{\circle*{.17}}}
\put(5,2){\red{\circle*{.17}}}
\put(6,2){\red{\circle*{.17}}}
\put(7,1){\red{\circle*{.17}}}
\put(8,1){\red{\circle*{.17}}}
\put(9,0){\red{\circle*{.17}}}
\end{picture}
\vspace*{3mm}
\caption{
   A Motzkin path of length 9,
   which gets weight $a_0^2 a_1 b_1^2 b_2 c_0 c_1 c_2$.
}
   \label{fig.motzkin}
\end{figure}

Now let ${\bf a} = (a_i)_{i \ge 0}$, ${\bf b} = (b_i)_{i \ge 1}$
and ${\bf c} = (c_i)_{i \ge 0}$ be indeterminates;
we will work in the ring $\Z[[{\bf a},{\bf b},{\bf c}]]$
of formal power series in these indeterminates.
We assign to each Motzkin path $\omega$
a weight $W(\omega) \in \Z[[{\bf a},{\bf b},{\bf c}]]$
that is the product of the weights for the individual steps,
where a rise starting at height~$i$ gets weight~$a_i$,
a fall starting at height~$i$ gets weight~$b_i$,
and a level step at height~$i$ gets weight~$c_i$
(see Figure~\ref{fig.motzkin}).

Define now for $k \ge 0$ the generating functions
\be
   f_k  \;=\;  \sum_{\omega \in \scrm_{k \to k}}  W(\omega)
   \;.
 \label{def.fk.flajolet}
\ee
These are well-defined elements of $\Z[[{\bf a},{\bf b},{\bf c}]]$
because there are finitely many $n$-step paths in $\scrm_{k \to k}$,
so each monomial occurs at most finitely many times.

Flajolet \cite{Flajolet_80} showed how to express
the generating functions $f_k$ as a continued fraction:

\begin{theorem}[Flajolet's master theorem]
   \label{thm.flajolet}
For each $k \ge 0$,
\be
   f_k
   \;=\;
   \cfrac{1}{1 - c_k - \cfrac{a_k b_{k+1}}{1 - c_{k+1} - \cfrac{a_{k+1} b_{k+2}}{1- \cdots}}}
 \label{eq.thm.flajolet}
\ee
as an identity in $\Z[[{\bf a},{\bf b},{\bf c}]]$.
\end{theorem}

\noindent
Of course, the identity \reff{eq.thm.flajolet} for one value of $k$
trivially implies it for all $k$, by redefining heights;
but in the proof it is natural to consider all $k$ simultaneously.

\par\bigskip\noindent{\sc Proof} \cite{Flajolet_80}.
Observe first that the right-hand side of \reff{eq.thm.flajolet}
is a well-defined element of $\Z[[{\bf a},{\bf b},{\bf c}]]$,
because all terms involving only
$(a_i)_{i \le k+r-1}$, $(b_i)_{i \le k+r}$ and $(c_i)_{i \le k+r-1}$
can be obtained by cutting off the continued fraction at level $r$,
yielding a rational fraction that expands into a well-defined
formal power series.

To prove \reff{eq.thm.flajolet}, we proceed as follows.
First define
\be
   f_k^\star  \;=\;  \sum_{\omega \in \scrm_{k \to k}^{\rm irred}}  W(\omega)
   \;,
\ee
where the sum is taken over {\em irreducible}\/ Motzkin paths at level~$k$,
i.e.\ paths of length $\ge 1$
that do not return to height~$k$ until the final step.
Since a Motzkin path can be uniquely decomposed as a concatenation of
some number $m \ge 0$ of irreducible Motzkin paths,
we have
\be
   f_k  \;=\;  \sum_{m=0}^\infty (f_k^\star)^m  \;=\;  {1 \over 1 - f_k^\star}
   \;.
 \label{eq.master.1}
\ee
On the other hand, an irreducible Motzkin path at level~$k$
is either a single level step at height~$k$ or else
begins with a rise $k \to k+1$ and ends with a fall $k+1 \to k$,
with an arbitrary Motzkin path at level~$k+1$ in-between; thus
\be
   f_k^\star \;=\;  c_k \,+\, a_k b_{k+1} f_{k+1}  \;.
 \label{eq.master.2}
\ee
Putting together \reff{eq.master.1} and \reff{eq.master.2},
we have
\be
   f_k  \;=\;  {1 \over 1 - c_k - a_k b_{k+1} f_{k+1}}
   \;.
 \label{eq.master.3}
\ee
Iterating \reff{eq.master.3}, we obtain \reff{eq.thm.flajolet}.
\qed

%

\bigskip

Let us now generalize this setup slightly by defining,
for any $k,\ell \ge 0$,
a \textbfit{Motzkin path at level $\bm{k \to \ell}$}
to be a path in $\Z \times \N$,
starting at height~$k$ and ending at height~$\ell$,
that stays always at height $\ge \min(k,\ell)$.
We write $\scrm_{k \to \ell}$
for the set of all Motzkin paths at level $k \to \ell$
that start at $(0,k)$.
For $\ell=k$ this reduces to the previous definition.
We then define the generating function
\be
   g_{k \to \ell}  \;=\;  \sum_{\omega \in \scrm_{k \to \ell}}  W(\omega)
   \;.
 \label{def.gkl.flajolet}
\ee

The generating functions $g_{k \to \ell}$ have a simple expression
in terms of the $f_k$:

\begin{proposition}
   \label{prop.gkl}
For $k,\ell \ge 0$ we have
\be
   g_{k \to \ell}
   \;=\;
   \begin{cases}
       f_k a_k f_{k+1} a_{k+1} \cdots f_{\ell-1} a_{\ell-1} f_\ell
                 & \textrm{if $k \le \ell$} \\[1mm]
       f_k b_k f_{k-1} b_{k-1} \cdots f_{\ell+1} b_{\ell+1} f_\ell
                 & \textrm{if $k \ge \ell$}
   \end{cases}
 \label{eq.prop.gkl}
\ee
\end{proposition}

\par\bigskip\noindent{\sc Proof}
\cite[pp.~295--296]{Goulden_83} \cite[pp.~II-7--II-8]{Viennot_83}.
For $k < \ell$, any path in $\scrm_{k \to\ell}$
can be uniquely decomposed by cutting it at its last return to height~$k$,
then at its last return to height~$k+1$, \ldots,
and so forth through its last return to height~$\ell-1$.
The pieces of this decomposition are an arbitrary Motzkin path
at level~$k$,
followed by a rise $k \to k+1$,
followed by an arbitrary Motzkin path at level~$k+1$,
followed by a rise $k+1 \to k+2$, \ldots,
followed by an arbitrary Motzkin path at level~$\ell$.

A similar argument handles the case $k > \ell$.
\qed

We can now specialize the foregoing results
to interpret continued fractions of the general form \reff{def.cfgen}.
Indeed, by taking $a_i = 1$, $b_i = A_i(t)$ and $c_i = \Delta_{i+1}(t)$,
we see that \reff{def.cfgen} is $A_0(t)$ times
the generating function for Motzkin paths at level~0
with the above weights.
Furthermore, the recurrence \reff{eq.master.3} relating $f_k$ to $f_{k+1}$
is identical to the recurrence \reff{eq.recurrence.fk.cfgen};
so the series $(f_k)_{k \ge 0}$ arising in our algorithm
are identical to those defined in \reff{def.fk.flajolet},
which enumerate Motzkin paths at level~$k$.
And finally, by Proposition~\ref{prop.gkl},
the series $(g_k/g_{-1})_{k \ge 0}$ arising in our refined algorithm
are identical to $(g_{0 \to k})_{k \ge 0}$  defined in \reff{def.gkl.flajolet},
which enumerate Motzkin paths at level~$0 \to k$.
We can therefore state:

\begin{proposition}
   \label{prop.cfgen}
The continued fraction \reff{def.cfgen} is $A_0(t)$ times
the generating function for Motzkin paths at level~$0$
in which each rise gets weight~$1$,
each fall starting at height~$i$ gets weight~$A_i(t)$,
and each level step at height~$i$ gets weight~$\Delta_{i+1}(t)$.

Moreover, $f_k$ is the generating function for
Motzkin paths at level~$k$ with these weights,
and $g_k$ ($k \ge 0$) is $g_{-1}(t)$ times the generating function for
Motzkin paths at level~$0 \to k$ with these weights.
\end{proposition}

Specializing this result we obtain interpretations of \reff{def.Jtype}
and \reff{def.Stype}.
In the latter case the level steps get weight $c_i = 0$,
so the relevant paths are Dyck paths.

\bigskip

Theorem~\ref{thm.flajolet} provides a powerful tool for proving
continued fractions in enumerative combinatorics.
Suppose that $P_n({\bf x})$ is the generating polynomial
for some class~$\scro_n$ of combinatorial objects of ``size $n$''
with respect to some set of statistics.
({\sc Example:} The polynomials $P_n(a,b)$ defined in \reff{eq.Pnab},
 which enumerate the set $\Sym_n$ of permutations of $\{1,\ldots,n\}$
 with respect to records and exclusive antirecords.)
And suppose that we can find a bijection from $\scro_n$
to some set $\scrl_n$ of {\em labeled}\/ Motzkin paths,
i.e.\ Motzkin paths augmented by putting labels on the steps,
where the label for a rise (resp.~fall, level step)
starting at height~$i$
belongs to some specified set $\scra_i$ (resp.~$\scrb_i$, $\scrc_i$)
of allowed labels.
Then the weights $a_i,b_i,c_i$ in the continued fraction \reff{eq.thm.flajolet}
can be obtained by summing over the labels.
This method goes back to Flajolet \cite{Flajolet_80};
for a detailed presentation with application to
permutations and set partitions,
see \cite[Sections~5--7]{Sokal-Zeng_masterpoly}.

\section{Connection with the work of Stieltjes and Rogers}

{}From now on we restrict attention to regular C-fractions
\be
   \cfrac{1}{1 - \cfrac{\alpha_1 t}{1 - \cfrac{\alpha_2 t}{1- \cdots}}}
 \label{eq.Sfrac.0}
\ee
and associated continued fractions
\be
   \cfrac{1}{1 - \gamma_0 t - \cfrac{\beta_1 t^2}{1 - \gamma_1 t - \cfrac{\beta_2 t^2}{1 - \cdots}}}
 \label{eq.Jfrac.0}
\ee
--- what combinatorialists call
\textbfit{S-fractions} and \textbfit{J-fractions}, respectively.

It is instructive to treat the coefficients $\balpha,\bbeta,\bgamma$
in these continued fractions as algebraic indeterminates.
We therefore write the S-fraction as
\be
   \cfrac{1}{1 - \cfrac{\alpha_1 t}{1 - \cfrac{\alpha_2 t}{1- \cdots}}}
   \;\;=\;\;
   \sum_{n=0}^\infty S_n(\balpha) \, t^n
 \label{eq.Stype.cfrac}
\ee
where $S_n(\balpha)$ is obviously a homogeneous polynomial of degree $n$
with nonnegative integer coefficients;
following Flajolet \cite{Flajolet_80},
we call it the \textbfit{Stieltjes--Rogers polynomial} of order $n$.
Likewise, we write the J-fraction as
\be
   \cfrac{1}{1 - \gamma_0 t - \cfrac{\beta_1 t^2}{1 - \gamma_1 t - \cfrac{\beta_2 t^2}{1 - \cdots}}}
   \;\;=\;\;
   \sum_{n=0}^\infty J_n(\bbeta,\bgamma) \, t^n
 \label{eq.Jtype.cfrac}
\ee
where $J_n(\bbeta,\bgamma)$
is a polynomial with nonnegative integer coefficients
that is quasi-homogeneous of degree $n$
if we assign weight~1 to each~$\gamma_i$ and weight~2 to each~$\beta_i$;
again following Flajolet \cite{Flajolet_80},
we call it the \textbfit{Jacobi--Rogers polynomial} of order $n$.

Since these are polynomials with nonnegative integer coefficients,
it is natural to ask what they count.
Flajolet's master theorem provides the immediate answer:

\begin{theorem}[Combinatorial interpretation of J-fractions and S-fractions]
\hfill\break
\vspace*{-7mm}
\begin{itemize}
   \item[(a)]  The Jacobi--Rogers polynomial $J_n(\bbeta,\bgamma)$
is the generating polynomial for Motzkin paths of length $n$,
in which each rise gets weight~1,
each fall from height~$i$ gets weight $\beta_i$,
and each level step at height~$i$ gets weight $\gamma_i$.
   \item[(b)]  The Stieltjes--Rogers polynomial $S_n(\balpha)$
is the generating polynomial for Dyck paths of length $2n$,
in which each rise gets weight~1
and each fall from height~$i$ gets weight $\alpha_i$.
\end{itemize}
\end{theorem}

\noindent
(We here made the arbitrary choice to weight the falls and not the rises.
Of course we could have done the reverse.)

But we can go farther.
Let us define a \textbfit{partial Motzkin path}
to be a path in the upper half-plane $\Z \times \N$,
starting on the horizontal axis but ending anywhere,
using the steps $(1,1)$, $(1,0)$ and $(1,-1)$.
Now define the \textbfit{generalized Jacobi--Rogers polynomial}
$J_{n,k}(\bbeta,\bgamma)$ to be the generating polynomial
for partial Motzkin paths from $(0,0)$ to $(n,k)$,
in which each rise gets weight~1,
each fall from height~$i$ gets weight $\beta_i$,
and each level step at height~$i$ gets weight $\gamma_i$.
Obviously $J_{n,k}$ is nonvanishing only for $0 \le k \le n$,
so we have an infinite lower-triangular array
$\sfJ = \big( J_{n,k}(\bbeta,\bgamma) \big)_{\! n,k \ge 0}$
in which the zeroth column displays
the ordinary Jacobi--Rogers polynomials $J_{n,0} = J_n$.
On the diagonal we have $J_{n,n} = 1$,
and on the first subdiagonal we have $J_{n,n-1} = \sum_{i=0}^{n-1} \gamma_i$.
By considering the last step of the path, we see that
the polynomials $J_{n,k}(\bbeta,\bgamma)$ satisfy the recurrence
\be
   J_{n+1,k}
   \;=\;
   J_{n,k-1}  \:+\:  \gamma_k J_{n,k} \:+\:  \beta_{k+1} J_{n,k+1}
 \label{eq.Jnk.recursion}
\ee
with the initial condition $J_{0,k} = \delta_{k0}$
(where of course we set $J_{n,-1} = 0$).

Similarly, let us define a \textbfit{partial Dyck path}
to be a partial Motzkin path without level steps.
Clearly a partial Dyck path starting at the origin
must stay on the even sublattice.
Now define the
\textbfit{generalized Stieltjes--Rogers polynomial of the first kind}
$S_{n,k}(\balpha)$
to be the generating polynomial for Dyck paths
starting at $(0,0)$ and ending at $(2n,2k)$,
in which each rise gets weight 1
and each fall from height~$i$ gets weight $\alpha_i$.
Obviously $S_{n,k}$ is nonvanishing only for $0 \le k \le n$,
so we have an infinite lower-triangular array
$\sfS = (S_{n,k}(\balpha))_{n,k \ge 0}$
in which the zeroth column displays
the ordinary Stieltjes--Rogers polynomials $S_{n,0} = S_n$.
We have $S_{n,n} = 1$ and $S_{n,n-1} = \sum_{i=1}^{2n-1} \alpha_i$.

Likewise, let us define
the \textbfit{generalized Stieltjes--Rogers polynomial of the second kind}
$S'_{n,k}(\balpha)$
to be the generating polynomial for Dyck paths
starting at $(0,0)$ and ending at $(2n+1,2k+1)$,
in which again each rise gets weight 1
and each fall from height~$i$ gets weight $\alpha_i$.
Since $S'_{n,k}$ is nonvanishing only for $0 \le k \le n$,
we obtain a second infinite lower-triangular array
$\sfS' = (S'_{n,k}(\balpha))_{n,k \ge 0}$.
We have $S'_{n,n} = 1$ and $S'_{n,n-1} = \sum_{i=1}^{2n} \alpha_i$.

The polynomials $S_{n,k}(\balpha)$ and $S'_{n,k}(\balpha)$
manifestly satisfy the joint recurrence
\begin{subeqnarray}
   S'_{n,k}   & = &  S_{n,k} \:+\: \alpha_{2k+2} \, S_{n,k+1}
       \\[2mm]
   S_{n+1,k}  & = &  S'_{n,k-1} \:+\: \alpha_{2k+1} \, S'_{n,k}
 \slabel{eq.SnlSprimenl.recurrence.b}
 \label{eq.SnlSprimenl.recurrence}
\end{subeqnarray}
for $n,k \ge 0$,
with the initial conditions $S_{0,k} = \delta_{k 0}$ and $S'_{n,-1} = 0$.
It follows that the $S_{n,k}$ satisfy the recurrence
\be
   S_{n+1,k}
   \;=\;
   S_{n,k-1}
      \:+\:  (\alpha_{2k} + \alpha_{2k+1}) \, S_{n,k}
      \:+\:  \alpha_{2k+1} \alpha_{2k+2} \, S_{n,k+1}
 \label{eq.Snl.recursion}
\ee
(where $S_{n,-1} = 0$ and $\alpha_0 = 0$),
while the $S'_{n,k}$ satisfy the recurrence
\be
   S'_{n+1,k}
   \;=\;
   S'_{n,k-1}
      \:+\:  (\alpha_{2k+1} + \alpha_{2k+2}) \, S'_{n,k}
      \:+\:  \alpha_{2k+2} \alpha_{2k+3} \, S'_{n,k+1}
   \;.
 \label{eq.Snlprime.recursion}
\ee
Note that \reff{eq.Snl.recursion} and \reff{eq.Snlprime.recursion}
have the same form as \reff{eq.Jnk.recursion},
when $\bbeta$ and $\bgamma$ are defined suitably in terms of the $\balpha$:
these correspondences are examples of \textbfit{contraction formulae}
\cite[p.~21]{Wall_48} \cite[p.~V-31]{Viennot_83}
that express an S-fraction as an equivalent J-fraction.
The recurrences
\reff{eq.Jnk.recursion}/\reff{eq.Snl.recursion}/\reff{eq.Snlprime.recursion}
define implicitly the (tridiagonal) \textbfit{production matrices}
for $\sfJ$, $\sfS$ and $\sfS'$:
see \cite{Deutsch_05,Deutsch_09,latpath_SRTR}.
Some workers call the arrays $\sfJ$, $\sfS$ and/or $\sfS'$
the \textbfit{Stieltjes table}.


The columns of the arrays $\sfS$ and $\sfS'$
are closely related to the series $g_k(t)$ of the
Euler--Viscovatov algorithm \reff{eq.alg.gk}/\reff{eq.alg.gkn}
with $g_{-1} = 1$
for the S-fraction \reff{eq.Sfrac.0}, as follows:

\begin{proposition}
   \label{prop.gk.Snk}
Let $g_{-1}(t) = 1$, let $g_0(t)$ be given by the S-fraction \reff{eq.Sfrac.0},
and let the series $(g_k)_{k \ge -1}$ satisfy the recurrence
\be
   g_k(t) - g_{k-1}(t)  \;=\; \alpha_{k+1} t \, g_{k+1}(t)
   \qquad\hbox{for $k \ge 0$}
   \;.
 \label{eq.recurrence.gk.Sfrac}
\ee
Then, in terms of the coefficients $g_{k,n}$ defined by
$g_k(t) = \sum\limits_{n=0}^\infty g_{k,n} t^n$,
we have
\begin{subeqnarray}
   g_{2j,n}    & = &  S_{n+j,j}
    \slabel{eq.prop.gk.Snk.even}   \\[2mm]
   g_{2j+1,n}  & = &  S'_{n+j,j}
  \label{eq.prop.gk.Snk}
\end{subeqnarray}
\end{proposition}

\noindent
In other words, the columns of $\sfS$ (resp.~$\sfS'$)
coincide with the coefficients of the even (resp.~odd) $g_k$,
but shifted downwards to start at the diagonal.

We will give two proofs of Proposition~\ref{prop.gk.Snk}:
one combinatorial and one algebraic.

\firstproof
We apply Flajolet's master theorem (Theorem~\ref{thm.flajolet})
with $a_i = 1$, $b_i = \alpha_i t$ and $c_i = 0$.
Then $f_0(t)$ is the S-fraction \reff{eq.Stype.cfrac},
and $f_k(t)$ is the analogous S-fraction but starting at $\alpha_{k+1}$.
By Proposition~\ref{prop.gkl} we have
$g_{0 \to \ell} = f_0 f_1 \cdots f_\ell$,
which equals the $g_\ell$ of the Euler--Gauss recurrence
\reff{eq.recurrence.gk.Sfrac} (since $g_{-1} = 1$).
So $g_\ell$ is the generating function for Dyck paths at level $0 \to \ell$
with the weights given above.
The coefficient of $t^n$ in $g_\ell$
corresponds to paths with $n$ falls and $n+\ell$ rises,
so the endpoint is $(2n+\ell,\ell)$.
If $\ell = 2j$, this gives $S_{n+j,j}$;
if $\ell = 2j+1$ this gives $S'_{n+j,j}$.
\qed

\secondproof
The recurrence \reff{eq.recurrence.gk.Sfrac} can be written
in terms of the coefficients $g_{k,n}$ as
\be
   g_{k,n} - g_{k-1,n}  \;=\;  \alpha_{k+1} \, g_{k+1,n-1}
   \;.
\ee
Evaluating this for $k=2j$ and $k=2j+1$ and using \reff{eq.prop.gk.Snk},
we recover the recurrences \reff{eq.SnlSprimenl.recurrence}.
Note also that $S'_{n,-1} = g_{-1,n+1} = 0$ by hypothesis.
\qed

\bigskip

\begin{example}
Consider the continued fraction \reff{eq.nfact.contfrac}.
With $g_{-1} = 1$, the first few $g_k$ are
\begin{subeqnarray}
g_0(t)  & = &  1+t+2 t^2+6 t^3+24 t^4+120 t^5+720 t^6+\ldots
   \\[1mm]
g_1(t)  & = &  1+2 t+6 t^2+24 t^3+120 t^4+720 t^5+5040 t^6+\ldots
   \\[1mm]
g_2(t)  & = &  1+4 t+18 t^2+96 t^3+600 t^4+4320 t^5+35280 t^6+\ldots
   \\[1mm]
g_3(t)  & = &  1+6 t+36 t^2+240 t^3+1800 t^4+15120 t^5+141120 t^6+\ldots
   \\[1mm]
g_4(t)  & = &  1+9 t+72 t^2+600 t^3+5400 t^4+52920 t^5+564480 t^6+\ldots
   \\[1mm]
g_5(t)  & = &  1+12 t+120 t^2+1200 t^3+12600 t^4+141120 t^5+1693440 t^6+\ldots
   \\[1mm]
g_6(t)  & = &  1+16 t+200 t^2+2400 t^3+29400 t^4+376320 t^5+5080320 t^6+\ldots
   \qquad\quad
\end{subeqnarray}
while the first few rows of $\sfS$ and $\sfS'$ are
\begin{eqnarray}
   \sfS
   & = &
   \begin{bmatrix}
 1 &   &   &   &   &   &      \\
 1 & 1 &   &   &   &   &      \\
 2 & 4 & 1 &   &   &   &      \\
 6 & 18 & 9 & 1 &   &   &      \\
 24 & 96 & 72 & 16 & 1 &   &      \\
 120 & 600 & 600 & 200 & 25 & 1 &      \\
 720 & 4320 & 5400 & 2400 & 450 & 36 & 1    \\
 \vdots & \vdots & \vdots & \vdots & \vdots & \vdots & \vdots
   \end{bmatrix}
       \\[5mm]
   \sfS'
   & = &
   \begin{bmatrix}
 1 &   &   &   &   &   &      \\
 2 & 1 &   &   &   &   &      \\
 6 & 6 & 1 &   &   &   &      \\
 24 & 36 & 12 & 1 &   &   &      \\
 120 & 240 & 120 & 20 & 1 &   &      \\
 720 & 1800 & 1200 & 300 & 30 & 1 &      \\
 5040 & 15120 & 12600 & 4200 & 630 & 42 & 1    \\
 \vdots & \vdots & \vdots & \vdots & \vdots & \vdots & \vdots
   \end{bmatrix}
\end{eqnarray}
The correspondences \reff{eq.prop.gk.Snk} can be observed.
{}From \reff{eq.euler.gk} we have
\begin{subeqnarray}
   g_{2j-1,n}  & = & \binom{n+j}{n} \binom{n+j-1}{n} \, n!
       \\[3mm]
   g_{2j,n}    & = & \binom{n+j}{n}^{\! 2} \, n!
 \label{eq.euler.gk.coeffs}
\end{subeqnarray}
and hence
\begin{subeqnarray}
   S_{n,k}  & = &  g_{2k,n-k}  \;=\;  \binom{n}{k}^{\! 2} \, (n-k)!
       \\[2mm]
   S'_{n,k}  & = &  g_{2k+1,n-k}  \;=\;  \binom{n+1}{k+1} \, \binom{n}{k} \, (n-k)!
\end{subeqnarray}
The recurrences \reff{eq.SnlSprimenl.recurrence}--\reff{eq.Snlprime.recursion}
with $\alpha_{2j-1} = \alpha_{2j} = j$ can easily be checked.
\myendremark
\end{example}

\medskip

{\bf Exercise.}  Work out the corresponding formulae for the
continued fractions \reff{eq.xupperfact.contfrac} and \reff{def.contfrac.Bxy}.
\myendremark

\bigskip

An analogous result connects the series $g_k(t)$ of the
Euler--Viscovatov algorithm \reff{eq.alg.gk.Jtype}
with $g_{-1} = 1$
for the J-fraction \reff{eq.Jfrac.0}
to the columns of the matrix $\sfJ$:

\begin{proposition}
   \label{prop.gk.Jnk}
Let $g_{-1}(t) = 1$, let $g_0(t)$ be given by the J-fraction \reff{eq.Jfrac.0},
and let the series $(g_k)_{k \ge -1}$ satisfy the recurrence
\be
   g_k(t) - g_{k-1}(t)  \;=\;
     \gamma_k t \, g_k(t)  \:+\:  \beta_{k+1} t^2 \, g_{k+1}(t)
   \qquad\hbox{for $k \ge 0$}
   \;.
 \label{eq.recurrence.gk.Jfrac}
\ee
Then, in terms of the coefficients $g_{k,n}$ defined by
$g_k(t) = \sum\limits_{n=0}^\infty g_{k,n} t^n$,
we have
\be
   g_{k,n}  \;=\;  J_{n+k,k}
   \;.
  \label{eq.prop.gk.Jnk}
\ee
\end{proposition}

Once again, this can be proven either combinatorially or algebraically;
these are left as exercises for the reader.

\bigdash

We can also interpret the
{\em exponential}\/ generating functions
of the columns of these lower-triangular arrays,
by using Hankel matrices.
Given a sequence $\ba = (a_n)_{n \ge 0}$ and an integer $m \ge 0$,
we define the $m$-shifted infinite Hankel matrix
$H_\infty^{(m)}(\ba) = (a_{i+j+m})_{i,j \ge 0}$.
We will apply this to the sequences
$\bJ = (J_n(\bbeta,\bgamma))_{n \ge 0}$
and $\bS = (S_n(\balpha))_{n \ge 0}$
of Jacobi--Rogers and Stieltjes--Rogers polynomials.
It turns out that the corresponding Hankel matrices
have beautiful $LDL^{\rm T}$ factorizations
in terms of the triangular arrays
of generalized Jacobi--Rogers and Stieltjes--Rogers polynomials:

\begin{theorem}[$LDL^{\rm T}$ factorization of Hankel matrices of Jacobi--Rogers and {Stieltjes}--Rogers polynomials]
   \label{thm.factorization}
We have the factorizations
\begin{itemize}
   \item[(a)]  $H_\infty^{(0)}(\bJ) \,=\,  \sfJ D \sfJ^{\rm T}$ where
$D = \diag(1, \beta_1, \beta_1 \beta_2, \ldots)$;
   \item[(b)]  $H_\infty^{(0)}(\bS) \,=\,  \sfS D \sfS^{\rm T}$ where
$D = \diag(1, \alpha_1 \alpha_2, \alpha_1 \alpha_2 \alpha_3 \alpha_4, \ldots)$;
   \item[(c)]  $H_\infty^{(1)}(\bS) \,=\,  \sfS' D' (\sfS')^{\rm T}$ where
$D = \diag(\alpha_1, \alpha_1 \alpha_2 \alpha_3,
                     \alpha_1 \alpha_2 \alpha_3 \alpha_4 \alpha_5, \ldots)$.
\end{itemize}
\end{theorem}

\proof
It suffices to note the identity
\cite[p.~351]{Aigner_01a} \cite[Remark~2.2]{Ismail_10}
\be
   J_{n+n',0}(\bbeta,\bgamma)
   \;=\;
   \sum_{k=0}^\infty
   J_{n,k}(\bbeta,\bgamma)
   \biggl( \prod_{i=1}^k \beta_i \! \biggr)
   J_{n',k}(\bbeta,\bgamma)
   \;,
 \label{eq.Jnl.hankel.0}
\ee
which arises from splitting a Motzkin path of length $n+n'$
into its first $n$ steps and its last $n'$ steps,
and then imagining the second part run backwards:
the factor $\prod_{i=1}^k \beta_i$ arises from the fact that
when we reversed the path we interchanged rises with falls
and thus lost a factor $\prod_{i=1}^k \beta_i$
for those falls that were not paired with rises.
The identity \reff{eq.Jnl.hankel.0} can be written in matrix form as
in part~(a).

The proofs of (b) and (c) are similar.
\qed

We can now prove an important equivalent formulation
of the factorization $H_\infty^{(0)}(\bJ) = \sfJ D \sfJ^{\rm T}$,
known as \textbfit{Rogers' addition formula} \cite{Rogers_07}.
We start with a simple observation:

\begin{lemma}[Bivariate egf of a Hankel matrix]
   \label{lemma.egf.hankel}
Let $\ba = (a_n)_{n \ge 0}$ be a sequence in
a commutative ring $R$ containing the rationals,
and let
\be
   A(t)  \;=\;  \sum_{n=0}^\infty a_n \, {t^n \over n!}
\ee
be its exponential generating function.
Then
\be
   A(t+u)
   \;=\;
   \sum_{n,n'=0}^\infty a_{n+n'} \: {t^n \over n!} \, {u^{n'} \over n'!}
   \;.
\ee
That is, $A(t+u)$ is the bivariate exponential generating function
of the Hankel matrix $H_\infty^{(0)}(\ba)$.
\end{lemma}

\proof
An easy computation.
\qed

As an immediate consequence, we get:

\begin{corollary}
   \label{cor.Jtype.converse.rogers}
Let $L = (\ell_{nk})_{n,k \ge 0}$ be a lower-triangular matrix
with entries in a commutative ring $R$ containing the rationals,
let
\be
   L_k(t)  \;=\;  \sum_{n=k}^\infty \ell_{nk} \, {t^n \over n!}
\ee
be the exponential generating function of the $k$th column of $L$,
and let $D = \diag(d_0,d_1,\ldots)$ be a diagonal matrix with entries in~$R$.
Let $\ba = (a_n)_{n \ge 0}$ be a sequence in $R$, and let
\be
   A(t)  \;=\;  \sum_{n=0}^\infty a_n \, {t^n \over n!}
\ee
be its exponential generating function.
Then $LDL^{\rm T} = H_\infty^{(0)}(\ba)$ if and only if
\be
   A(t+u)
   \;=\;
   \sum\limits_{k=0}^\infty d_k \, L_k(t) \, L_k(u)
   \;.
 \label{eq.cor.Jtype.converse.rogers}
\ee
\end{corollary}


On the other hand, a converse to the factorization
of Theorem~\ref{thm.factorization}(a) can be proven.
We recall that an element of a commutative ring $R$ is called
\textbfit{regular} if it is neither zero nor a divisor of zero,
and that a diagonal matrix is called regular if all its diagonal elements are.
We then have the following result
(\!\!\cite{Sokal_totalpos}, based on
 \cite[Theorem~1]{Peart_00} \cite[Theorem~2.1]{Woan_01}),
which we state here without proof:

\begin{proposition}
   \label{prop.Jtype.converse.bis2}
Let $R$ be a commutative ring,
let $L$ be a unit-lower-triangular matrix with entries in $R$,
let $D = \diag(d_0,d_1,\ldots)$ be a {\em regular} diagonal matrix
with entries in $R$,
and let $\ba = (a_n)_{n \ge 0}$ be a sequence in $R$.
If $L D L^{\rm T} = H_\infty^{(0)}(\ba)$,
then there exist sequences $\bbeta = (\beta_n)_{n \ge 1}$ and
$\bgamma = (\gamma_n)_{n \ge 0}$ in $R$
such that $d_n = d_0 \beta_1 \cdots \beta_n$,
$L = \sfJ(\bbeta,\bgamma)$ and $\ba = d_0 \bJ(\bbeta,\bgamma)$.
In particular, $\ba$ equals $d_0$ times the zeroth column of $L$.
\end{proposition}

Putting together Theorem~\ref{thm.factorization},
Corollary~\ref{cor.Jtype.converse.rogers}
and Proposition~\ref{prop.Jtype.converse.bis2},
we conclude (compare \cite[Theorem~53.1]{Wall_48}):

\begin{theorem}[Rogers' addition formula]
   \label{thm.rogers}
The column exponential generating functions of
the matrix of generalized Jacobi--Rogers polynomials,
\be
   \scrj_k(t;\bbeta,\bgamma)
   \;\eqdef\;
   \sum_{n=k}^\infty J_{n,k}(\bbeta,\bgamma) \, {t^n \over n!}
   \;,
\ee
satisfy
\be
   \scrj_0(t+u;\bbeta,\bgamma)
   \;=\;
   \sum\limits_{k=0}^\infty \beta_1 \cdots \beta_k \,
           \scrj_k(t;\bbeta,\bgamma) \, \scrj_k(u;\bbeta,\bgamma)
   \;.
\ee

And conversely, if $A(t)$ and $F_0(t), F_1(t), \ldots$ are formal power series
(with elements in a commutative ring $R$ containing the rationals)
satisfying
\be
   A(t) \;=\; 1 + O(t) \,,\qquad
   F_k(t) \;=\; {t^k \over k!} \,+\, \mu_k {t^{k+1} \over (k+1)!} \,+\,
                  O(t^{k+2})
 \label{eq.thm.rogers.1}
\ee
and
\be
   A(t+u)
   \;=\;
   \sum\limits_{k=0}^\infty \beta_1 \cdots \beta_k \, F_k(t) \, F_k(u)
 \label{eq.thm.rogers.2}
\ee
for some {\em regular} elements $\bbeta = (\beta_k)_{k \ge 1}$,
then $A(t) = F_0(t)$ and $F_k(t) = \scrj_k(t;\bbeta,\bgamma)$
with the given $\bbeta$
and with $\gamma_k = \mu_k - \mu_{k-1}$ (where $\mu_{-1} \eqdef 0$).
\end{theorem}

Here the formula for $\gamma_k$ follows from
$J_{k+1,k} = \sum_{i=0}^{k} \gamma_i$.

\bigskip

\begin{example}
The \textbfit{secant numbers}\footnote{
   See \cite{Sokal_eulernumbers} and the references cited therein
   for more information concerning the secant numbers
   and the closely-related tangent numbers.
}
$E_{2n}$ are defined by the
exponential generating function
\be
   \sec t  \;=\;  \sum_{n=0}^\infty E_{2n} \: {t^{2n} \over (2n)!}
   \;.
\ee
More generally, the \textbfit{secant power polynomials} $E_{2n}(x)$
are defined by the exponential generating function
\be
   (\sec t )^x \;=\;  \sum_{n=0}^\infty E_{2n}(x) \: {t^{2n} \over (2n)!}
   \;.
\ee
{}From the high-school angle-addition formula
\begin{subeqnarray}
   \cos(t+u)
   & = &
   (\cos t)(\cos u) \,-\, (\sin t)(\sin u)
       \\[1mm]
   & = &
   (\cos t)(\cos u) \, [ 1 \,-\, (\tan t)(\tan u) ]
\end{subeqnarray}
we obtain
\be
   [\sec(t+u)]^x
   \;=\;
   (\sec t)^x (\sec u)^x \:
      \sum_{k=0}^\infty \binom{x+k-1}{k} \, (\tan t)^k \, (\tan u)^k
   \;,
\ee
which is of the form \reff{eq.thm.rogers.1}/\reff{eq.thm.rogers.2} with
\be
   \beta_k \:=\; k(x+k-1)
   \,,\qquad
   F_k(t) \:=\: {(\sec t)^x \, (\tan t)^k \over k!}
          \:=\: {t^k \over k!} \,+\, O(t^{k+2})
   \;,
\ee
so that $\mu_k = 0$ and hence $\gamma_k = 0$.
Theorem~\ref{thm.rogers} then implies that the ordinary generating function
of the secant power polynomials is given by the J-fraction
\be
   \sum_{n=0}^\infty E_{2n}(x) \, t^{2n}
   \;=\;
   \cfrac{1}{1 - \cfrac{1 \cdot x t^2}{1 - \cfrac{2(x+1) t^2}{1 - \cfrac{3(x+2) t^2}{1 - \cdots}}}}
   \;.
\ee
After renaming $t^2 \to t$, this is actually an S-fraction
with coefficients $\alpha_n = n(x+n-1)$.
This example is due to Stieltjes \cite{Stieltjes_1889}
and Rogers \cite{Rogers_07}.
\myendremark
\end{example}

\begin{example}
Let us use Rogers' addition formula to give a second proof of
Euler's continued fraction \reff{eq.xupperfact.contfrac}
for the sequence of rising powers $(a^{\overline{n}})_{n \ge 0}$.
This sequence has the exponential generating function
\be
   \sum_{n=0}^\infty a^{\overline{n}} \: {t^n \over n!}
   \;=\;
   \sum_{n=0}^\infty \binom{a+n-1}{n} \, t^n
   \;=\;
   (1-t)^{-a}
   \;,
\ee
which satisfies the addition formula
\begin{subeqnarray}
   (1-t-u)^{-a}
   & = &
   (1-t)^{-a} \, (1-u)^{-a} \, \Bigl[ 1 \,-\, {tu \over (1-t)(1-u)} \Bigr]^{-a}
        \\[2mm]
   & = &
   (1-t)^{-a} \, (1-u)^{-a} \sum_{k=0}^\infty \binom{a+k-1}{k} \,
        \Bigl( {t \over 1-t} \Bigr)^{\! k} \,
        \Bigl( {u \over 1-u} \Bigr)^{\! k}
   \,.
   \qquad\qquad
\end{subeqnarray}
This expansion is of the form \reff{eq.thm.rogers.1}/\reff{eq.thm.rogers.2}
with $\beta_k = k(a+k-1)$ and
\be
   F_k(t) \;=\; (1-t)^{-a} \, \Bigl( {t \over 1-t} \Bigr)^{\! k} \, {1 \over k!}
          \;=\; {t^k \over k!} \,+\, (k+1)(k+a) {t^{k+1} \over (k+1)!} \,+\,
                O(t^{k+2})
   \;,
\ee
so that $\mu_k = (k+1)(k+a)$ and hence $\gamma_k = 2k+a$.
Moreover, the J-fraction \reff{eq.Jfrac.0}
with $\beta_k = k(a+k-1)$ and $\gamma_k = 2k+a$
is connected by the contraction formula
\cite[p.~21]{Wall_48} \cite[p.~V-31]{Viennot_83}
\begin{subeqnarray}
   \gamma_0  & = &  \alpha_1
       \slabel{eq.contraction_even.coeffs.a}   \\
   \gamma_n  & = &  \alpha_{2n} + \alpha_{2n+1}  \qquad\hbox{for $n \ge 1$}
       \slabel{eq.contraction_even.coeffs.b}   \\
   \beta_n  & = &  \alpha_{2n-1} \alpha_{2n}
       \slabel{eq.contraction_even.coeffs.c}
 \label{eq.contraction_even.coeffs}
\end{subeqnarray}
with the S-fraction having $\alpha_{2k-1} = a+k-1$ and $\alpha_{2k} = k$.
This completes the proof of \reff{eq.xupperfact.contfrac}.

We also see from this proof that the generalized Jacobi--Rogers polynomials
for the J-fraction with $\beta_k = k(a+k-1)$ and $\gamma_k = 2k+a$ are
\begin{subeqnarray}
   J_{n,k}  \;\eqdef\; \Bigl[ {t^n \over n!} \Bigr] \, F_k(t)
   & = &
   {n! \over k!} \: [t^n] \, (1-t)^{-a} \, \Bigl( {t \over 1-t} \Bigr)^{\! k}
        \\[2mm]
   & = &
   {n! \over k!} \: [t^{n-k}] \, (1-t)^{-(a+k)}
        \\[2mm]
   & = &
   {n! \over k!} \, \binom{a+n-1}{n-k}
        \\[2mm]
   & = &
   \binom{n}{k} \, (a+k)^{\overline{n-k}}
   \;.
 \label{eq.risingpower.Jnk}
\end{subeqnarray}
These also coincide with the generalized Stieltjes--Rogers polynomials
of the first kind $S_{n,k}$ for the corresponding S-fraction
\reff{eq.xupperfact.contfrac},
since the contraction formula \reff{eq.contraction_even.coeffs}
corresponds combinatorially \cite[p.~V-31]{Viennot_83}
to grouping pairs of steps of the Dyck path to create a Motzkin path
living on even heights.
Then \reff{eq.risingpower.Jnk} agrees with \reff{eq.euler.gk.BIS.2j}
in view of \reff{eq.prop.gk.Snk.even}.
\myendremark
\end{example}


See \cite[pp.~203--207]{Wall_48} \cite{Ismail_10}
for further discussion of Rogers' addition formula
and its applications to the derivation of continued fractions.

\bigskip

{\bf Historical remarks.}
The generalized Stieltjes--Rogers polynomials $S_{n,k}$ and $S'_{n,k}$
were introduced by Stieltjes \cite{Stieltjes_1889} in 1889
(his notation is $\alpha_{k,n}$ and $\beta_{k,n}$):
he defined them by the recurrences \reff{eq.SnlSprimenl.recurrence}.
He then proved the factorizations in Theorem~\ref{thm.factorization}(b,c)
by considering the quadratic forms associated to the symmetric matrices
$\sfS D \sfS^{\rm T}$ and $\sfS' D' (\sfS')^{\rm T}$:
he used the recurrence to prove that
the matrix $\sfS D \sfS^{\rm T}$ is Hankel,
i.e.\ is $H_\infty^{(0)}(\bb)$ for some sequence $\bb = (b_n)_{n \ge 0}$;
then, using the previously known formula
(for which he cited Frobenius and Stickelberger
\cite{Frobenius_1879,Frobenius_1881})
relating the coefficients $\balpha = (\alpha_n)_{n \ge 1}$
in an S-fraction to the Hankel determinants
of the power-series coefficients $\ba = (a_n)_{n \ge 0}$,
he concluded that $\ba = \bb$.
Stieltjes went on to use this matrix-decomposition method
to determine several explicit continued fractions related
to trigonometric functions and Jacobian elliptic functions.
See also the summary of this work
given in Stieltjes' 1894 memoir \cite[pp.~J.18--J.19]{Stieltjes_1894},
where the matrix factorizations are made explicit.

The reformulation of Stieltjes' factorization as an addition formula
is due to Rogers \cite{Rogers_07} in 1907.

The interpretation of $J_{n,k}$, $S_{n,k}$ and $S'_{n,k}$
in terms of partial Motzkin and Dyck paths is post-Flajolet folklore;
it goes back at least to \cite[Theorem~2.1 and Remark~2.2]{Ismail_10}.
\myendremark

%
%


\section{Timing tests}  \label{sec.timing}

How do the primitive algorithm \reff{eq.alg.fk}
and the refined algorithm \reff{eq.alg.gk}/\reff{eq.alg.gkn}
compare in computational efficiency?

\begin{figure}[p]
\centering
\includegraphics[width=0.85\textwidth]{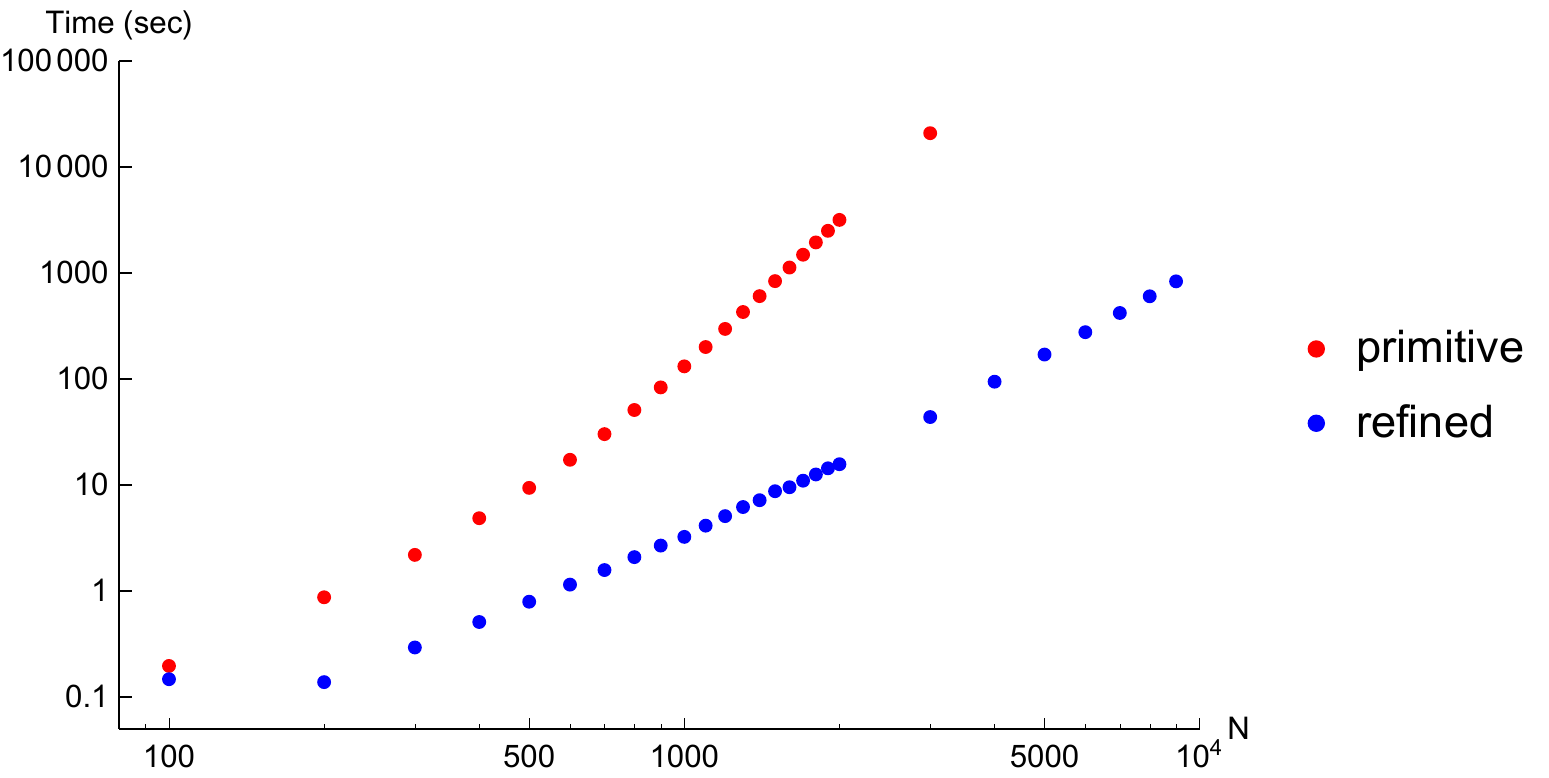}
\caption{
   Timings (in seconds) for the primitive algorithm (upper curve)
   and refined algorithm (lower curve)
   applied to the numerical series \reff{eq.nfact.contfrac}.
}
\label{fig.nfact}
\end{figure}

\begin{figure}[p]
\centering
\vspace*{-1cm}
\includegraphics[width=0.85\textwidth]{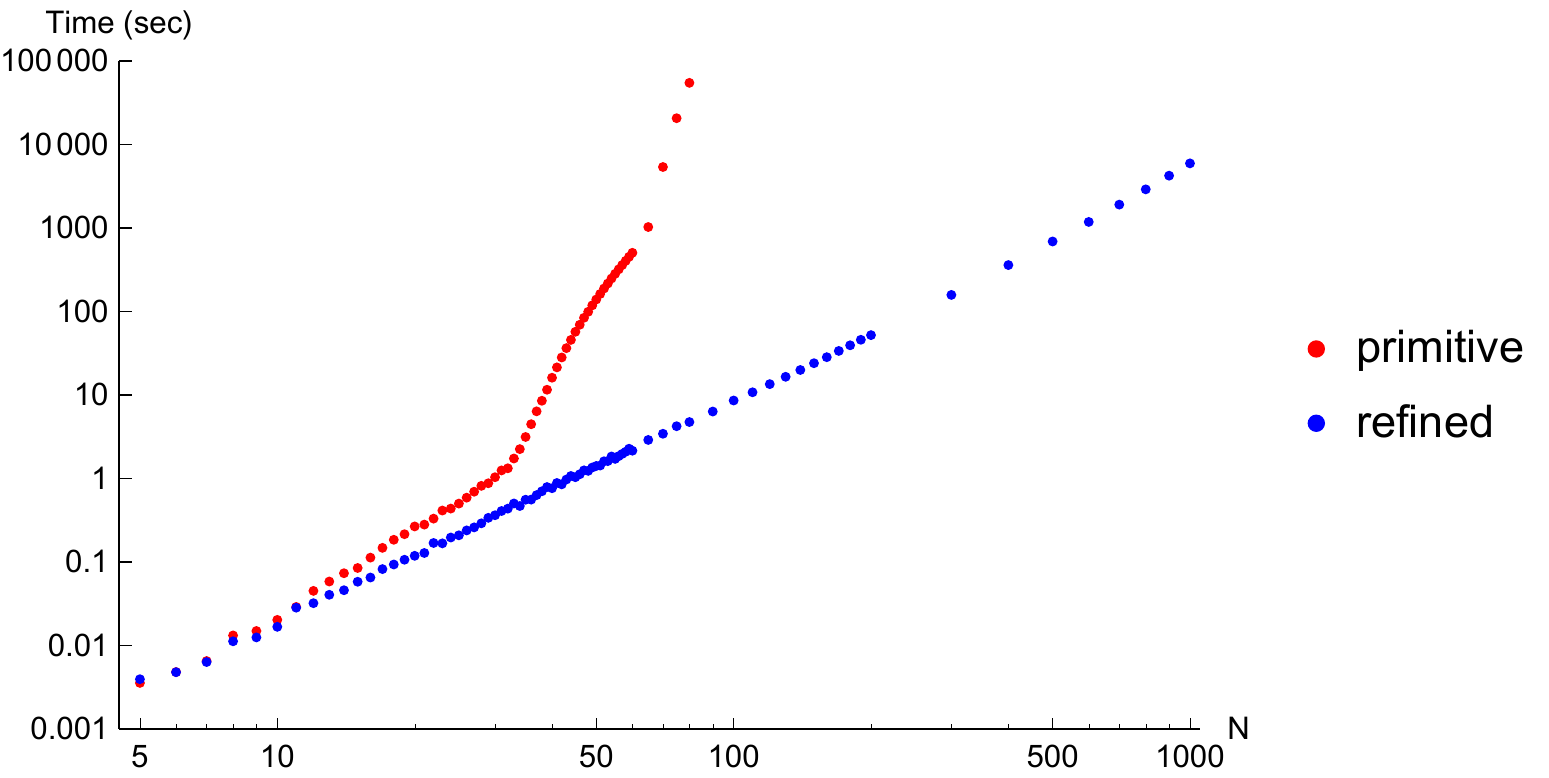} \\[2mm]
   (a) \\[2cm]
\includegraphics[width=0.85\textwidth]{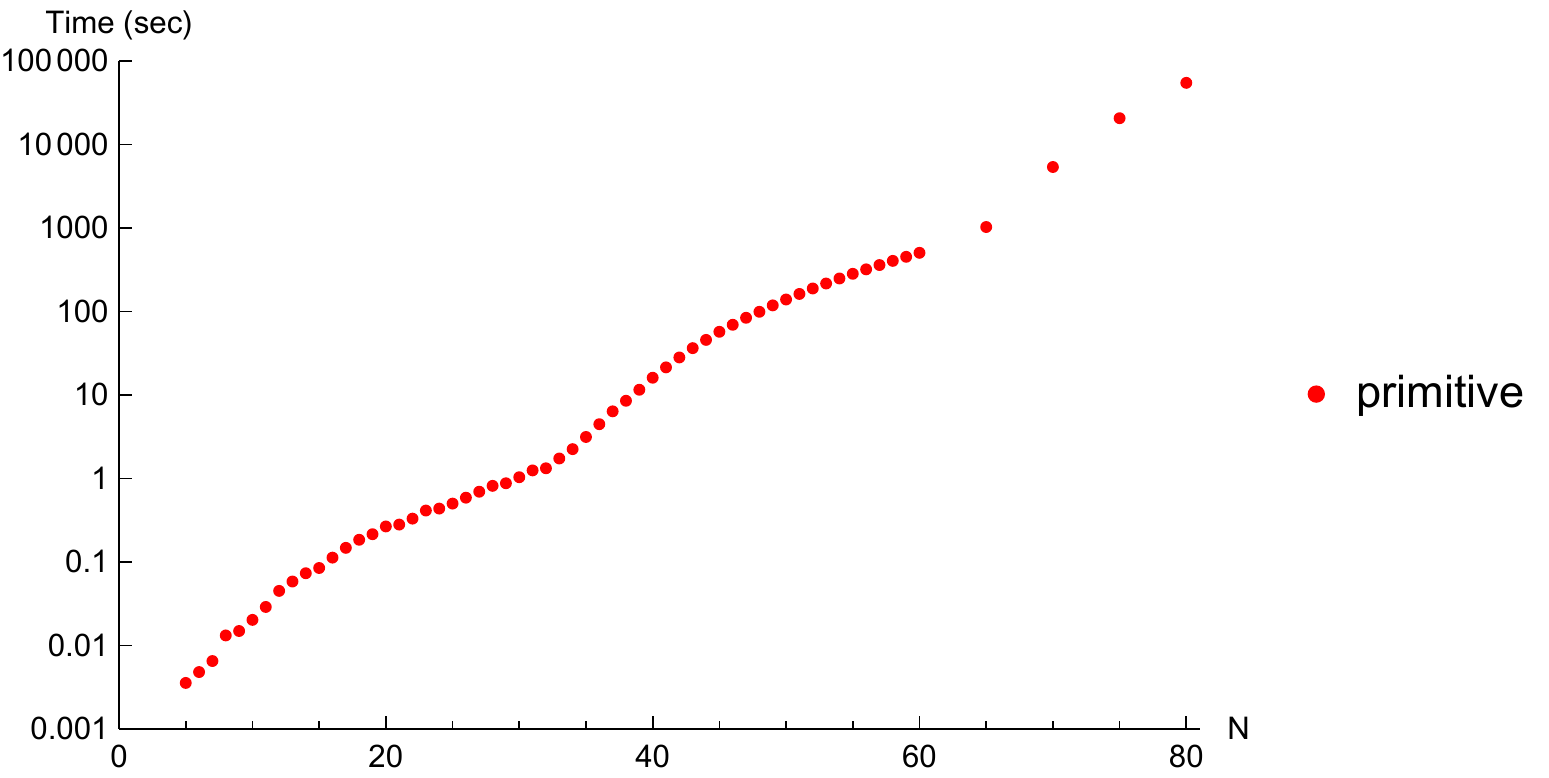} \\[2mm]
    (b) \\[5mm]
\caption{
   Timings (in seconds) for the primitive algorithm (upper curve)
   and refined algorithm (lower curve)
   applied to the polynomial series \reff{eq.xupperfact.contfrac}:
   log-log plot in (a), linear-log plot for the primitive algorithm in (b).
}
\label{fig.risingpower}
\end{figure}

Numerical timing experiments for the continued fractions
\reff{eq.nfact.contfrac} and \reff{eq.xupperfact.contfrac}
are reported in Table~\ref{table.nfact}/Figure~\ref{fig.nfact}
and Table~\ref{table.risingpower}/Figure~\ref{fig.risingpower},
respectively.
The computations were carried out in {\sc Mathematica} version~11.1
under Linux
on a machine with an Intel Xeon W-2133 CPU running at 3.60 GHz.
The primitive algorithm was programmed in both recursive and iterative forms;
the timings for the two versions were essentially identical.

For the numerical series $a_n = n!$,
the CPU time for the primitive algorithm behaves roughly like $N^{\approx 2}$
for the smaller values of $N$,
rising gradually to $N^{\approx 4.6}$ for $1000 \ltapprox N \ltapprox 3000$.
The CPU time for the refined algorithm behaves roughly like $N^{\approx 2}$
over the range $N \ltapprox 2000$,
rising gradually to $N^{\approx 2.8}$ for $N \approx 9000$.\footnote{
   Our computation for $N=10000$ required more memory than the available
   256 GB, which led to paging and an erratic timing;
   we have therefore suppressed this data point as unreliable.
}
This latter behavior is consistent with the theoretically expected
(but not yet reached) asymptotic CPU time of order $N^3 \log^2 N$,
arising as $\sim N^2$ field operations in \reff{eq.alg.gk}/\reff{eq.alg.gkn}
times a CPU time of order $N \log^2 N$ per operation:
here the operations are subtraction of numbers of magnitude roughly $N!$
(hence with $\sim N \log N$ digits)
and their division by the integers $\alpha_k$ of order $N$
(hence with $\sim \log N$ digits).
The advantage for the refined algorithm grows from a factor $\approx 6$
at $N=200$ to $\approx 500$ at $N=3000$.

For the polynomial series $a_n = a^{\overline{n}}$,
the CPU time for the primitive algorithm behaves roughly like $N^{\approx 3.3}$
for $5 \ltapprox N \ltapprox 30$,
bending suddenly at $N=30$ to a much more rapid growth $N^{\approx 10}$
[see Figure~\ref{fig.risingpower}(a)].
However, another possible interpretation is that the behavior
is exponential in $N$
[see Figure~\ref{fig.risingpower}(b)].
The CPU time for the refined algorithm, by contrast,
behaves like $N^{\approx 3}$ over the whole range $5 \le N \le 1000$,
with a slightly lower power ($\approx 2.7$) at the smallest values of $N$
and a slightly higher power ($\approx 3.1$) at the largest.
I~am not sure what should be the expected asymptotic behavior
for either algorithm.
The advantage for the refined algorithm grows from a factor $\approx 1.2$
at $N=10$ to $\approx 3$ at $N=30$ and $\approx 10000$ at $N=80$.

\begin{table}[p]
\centering
\footnotesize
\begin{tabular}{|r|r|r|r|}
\hline
              & \multicolumn{1}{|c|}{Primitive}
              & \multicolumn{1}{|c|}{Refined} &  \\
   \multicolumn{1}{|c|}{$N$} &
   \multicolumn{1}{|c|}{algorithm} &
   \multicolumn{1}{|c|}{algorithm} &
   \multicolumn{1}{|c|}{Ratio} \\
\hline
100  &  0.20  &  0.15  &  1.33    \\
200  &  0.87  &  0.14  &  6.32    \\
300  &  2.20  &  0.29  &  7.47    \\
400  &  4.87  &  0.51  &  9.53    \\
500  &  9.41  &  0.79  &  11.86    \\
600  &  17.32  &  1.15  &  15.06    \\
700  &  30.26  &  1.58  &  19.17    \\
800  &  51.10  &  2.09  &  24.44    \\
900  &  83.48  &  2.69  &  31.07    \\
1000  &  131.90  &  3.25  &  40.63    \\
1100  &  200.71  &  4.14  &  48.46    \\
1200  &  297.45  &  5.10  &  58.38    \\
1300  &  429.43  &  6.21  &  69.18    \\
1400  &  606.35  &  7.20  &  84.20    \\
1500  &  840.25  &  8.75  &  95.99    \\
1600  &  1128.79  &  9.54  &  118.28    \\
1700  &  1490.64  &  11.00  &  135.50    \\
1800  &  1947.84  &  12.59  &  154.68    \\
1900  &  2505.78  &  14.40  &  174.06    \\
2000  &  3176.93  &  15.74  &  201.85    \\
3000  &  20896.0\hphantom{0}  &  43.85  &  476.52    \\
4000  &        &  94.49   &              \\
5000  &        &  170.51   &              \\
6000  &        &  277.10   &              \\
7000  &        &  420.58   &              \\
8000  &        &  604.25   &              \\
9000  &        &  835.81   &              \\
\hline
\end{tabular}
\caption{
   Timings (in seconds) for the primitive and refined algorithms
   applied to the numerical series \reff{eq.nfact.contfrac}.
}
\label{table.nfact}
\end{table}

\begin{table}[p]
\centering
\vspace*{-3cm}
\begin{tabular}{|r|r|r|r|}
\hline
              & \multicolumn{1}{|c|}{Primitive}
              & \multicolumn{1}{|c|}{Refined} &  \\
   \multicolumn{1}{|c|}{$N$} &
   \multicolumn{1}{|c|}{algorithm} &
   \multicolumn{1}{|c|}{algorithm} &
   \multicolumn{1}{|c|}{Ratio} \\
\hline
10  &  0.02  &  0.02  &  1.21    \\
15  &  0.08  &  0.06  &  1.46    \\
20  &  0.27  &  0.12  &  2.25    \\
25  &  0.50  &  0.21  &  2.40    \\
30  &  1.04  &  0.36  &  2.85    \\
35  &  3.15  &  0.56  &  5.64    \\
40  &  16.13  &  0.77  &  21.07    \\
45  &  57.23  &  1.04  &  55.14    \\
50  &  139.52  &  1.41  &  98.66    \\
55  &  283.39  &  1.72  &  164.86    \\
60  &  505.61  &  2.15  &  234.67    \\
65  &  1029.79  &  2.90  &  355.29    \\
70  &  5390.53  &  3.44  &  1567.81    \\
75  &  20714.2\hphantom{0}  &  4.23  &  4893.62    \\
80  &  54919.5\hphantom{0}  &  4.75  &  11560.1\hphantom{0}    \\
90  &          &  6.35  &          \\
100  &          &  8.60  &          \\
110  &          &  10.79  &          \\
120  &          &  13.52  &          \\
130  &          &  16.54  &          \\
140  &          &  19.97  &          \\
150  &          &  24.06  &          \\
160  &          &  28.42  &          \\
170  &          &  33.76  &          \\
180  &          &  39.46  &          \\
190  &          &  45.91  &          \\
200  &          &  52.23  &          \\
300  &          &  158.25  &          \\
400  &          &  360.65  &          \\
500  &          &  691.27  &          \\
600  &          &  1184.81  &          \\
700  &          &  1910.57  &          \\
800  &          &  2909.85  &          \\
900  &          &  4244.91  &          \\
1000  &          &  5960.16  &          \\
\hline
\end{tabular}
\caption{
   Timings (in seconds) for the primitive and refined algorithms
   applied to the polynomial series \reff{eq.xupperfact.contfrac}.
}
\label{table.risingpower}
\end{table}

\bigskip

{\bf Some remarks.}
1.  When the primitive algorithm is programmed recursively
in {\sc Mathematica}, it is necessary to set {\tt \$Recursion Limit}
to a large enough number (or {\tt Infinity})
in order to avoid incomplete execution.

2.  Because of quirks in {\sc Mathematica}'s treatment of power series
with symbolic coefficients, the primitive algorithm (in either version)
applied to \reff{eq.xupperfact.contfrac} becomes exceedingly slow
for $N \gtapprox 10$ if the basic step is programmed simply as
\linebreak
{\tt f[k] = (1 - 1/f[k-1])/(alpha[k]*t)}.
Instead, it is necessary to write
\linebreak
{\tt f[k] = Map[Together, (1 - 1/f[k-1])/(alpha[k]*t)]}
in order to force the simplification
of rational-function expressions to polynomials.
I~thank Daniel Lichtblau for this crucial suggestion.
The results reported in Table~\ref{table.risingpower}
and Figure~\ref{fig.risingpower}
refer to this latter version of the program.

3. The timings reported here were obtained using {\sc Mathematica}'s
command {\tt Timing}, which under this operating system
apparently includes the total CPU time in all threads.
The real time elapsed was in some instances up to a factor $\approx 2$
smaller than this, due to partially parallel execution on this multi-core CPU.

4. One might wonder: Why on earth would one {\em want}\/ to compute
1000 or more continued-fraction coefficients?
One answer (perhaps not the only one)
is that the nonnegativity of the S-fraction coefficients $\alpha_n$
is a necessary and sufficient condition for a sequence $\ba = (a_n)_{n \ge 0}$
of real numbers to be a \textbfit{Stieltjes moment sequence},
i.e.~the moments of a positive measure on $[0,\infty)$;
this was shown by Stieltjes \cite{Stieltjes_1894} in 1894.
On the other hand, it is easy to concoct sequences
that are {\em not}\/ Stieltjes moment sequences
but which have $\alpha_n > 0$ until very high order.
Consider, for instance, the sequence\footnote{
   A closely related form of $a_n$ was suggested to me by
   Andrew Elvey Price \cite{Elvey_private}.
}
\be
   a_n  \;\eqdef\;  (1+\epsilon) \, n! \,-\, {\epsilon \over (n+1)^2}
   \;=\;  \int\limits_0^\infty x^n \:
             \Bigl[ (1+\epsilon) e^{-x} \,+\, \epsilon \log x \Bigr] \: dx
   \;,
\ee
which fails to be a Stieltjes moment sequence whenever $\epsilon > 0$
because the density is negative near $x=0$
(apply \cite[Corollary~2.10]{Sokal_eulernumbers}).
For $\epsilon = 1$, the first negative coefficient $\alpha_n$ is $n=6$;
for $\epsilon = 1/2$ it is $n = 20$;
for $\epsilon = 1/4$ it is $n = 178$;
for $\epsilon = 1/8$ it is some unknown (to me) $n > 1500$.
So it can be important to compute S-fraction coefficients to very high order
when trying to determine empirically whether a given sequence is or is not
a Stieltjes moment sequence.
\myendremark

\section{Final remarks}

The algorithm presented here is intended, in the first instance,
for use in exact arithmetic:  the field $F$ could be (for example)
the field $\Q$ of rational numbers, or more generally
the field $\Q(x_1,\ldots,x_n)$ of rational fractions in
indeterminates $x_1,\ldots,x_n$ with coefficients in $\Q$.
I~leave it to others to analyze the numerical (in)stability
of this algorithm when carried out in $F = \R$ or $\C$
with finite-precision arithmetic, and/or to devise alternative algorithms
with improved numerical stability.

The continued fractions discussed here are what could be called
{\em classical continued fractions}\/.
Very recently combinatorialists have developed a theory
of {\em branched continued fractions}\/,
based on generalizing Flajolet's master theorem (Theorem~\ref{thm.flajolet})
to other classes of lattice paths.
This idea was suggested by Viennot \cite[section~V.6]{Viennot_83},
carried forward in the Ph.D.~theses of
Roblet \cite{Roblet_94} and Varvak \cite{Varvak_04},
and then comprehensively developed by
P\'etr\'eolle, Sokal and Zhu \cite{latpath_SRTR,latpath_lah}.
There is a corresponding generalization of the Euler--Gauss recurrence method:
for instance, for the $m$-S-fractions,
which generalize the regular C-fractions,
the recurrence \reff{eq.recurrence.gk.Sfrac} is generalized to
\be
   g_k(t) - g_{k-1}(t)  \;=\; \alpha_{k+m} t \, g_{k+m}(t)
   \qquad\hbox{for } k \ge 0
 \label{eq.recurrence.gkm.0}
\ee
for a fixed integer $m \ge 1$.
Furthermore, Gauss' \cite{Gauss_1813} continued fraction
for the ratio of contiguous hypergeometric functions $\tfo$
can be generalized to $\FHyper{r}{s}$ for arbitrary $r,s$,
where now $m = \max(r-1,s)$;
the proof is based on \reff{eq.recurrence.gkm.0}.
See \cite{latpath_SRTR} for details on all of this,
and \cite{latpath_lah} for further applications.
On~the other hand, branched continued fractions are highly nonunique,
and I~do not know any algorithm for computing them.

\section*{Appendix}
\setcounter{section}{1}
\setcounter{equation}{0}
\def\theequation{\Alph{section}.\arabic{equation}}

Answer to the exercise posed in Section~\ref{sec.examples.2}:
\begin{subeqnarray}
   g_{2j-1}(t)
   & = &
   \sum_{n=0}^\infty \sum_{k=0}^n \stirlingsubset{n+j}{k+j} \, \binom{k+j-1}{k}
                                                 \, x^k y^{n-k} \, t^n
       \\[3mm]
   g_{2j}(t)
   & = &
   \sum_{n=0}^\infty \sum_{k=0}^n \stirlingsubset{n+j}{k+j} \, \binom{k+j}{k}
                                                 \, x^k y^{n-k} \, t^n
\end{subeqnarray}

\section*{Acknowledgments}
\addcontentsline{toc}{section}{Acknowledgments}

I wish to thank Gaurav Bhatnagar, Bishal Deb, Bill Jones, Xavier Viennot
and Jiang Zeng for helpful conversations and/or correspondence.
I~am especially grateful to Gaurav Bhatnagar for reemphasizing to me
the power and elegance of the Euler--Gauss recurrence method,
and for drawing my attention to Askey's masterful survey \cite{Askey_89}
as well as to his own wonderful survey article \cite{Bhatnagar_14}.
I~also thank Daniel Lichtblau for help with {\sc Mathematica}.

This research was supported in part by
the U.K.~Engineering and Physical Sciences Research Council grant EP/N025636/1.

\end{document}